\documentclass{autart}    
\usepackage{amsmath}
\usepackage{amssymb} 
\usepackage{tikz}

\usepackage{graphicx}          

\begin{document}

\begin{frontmatter}

\title{Disturbance-to-state stabilization by output feedback of nonlinear ODE cascaded with a reaction-diffusion equation} 

\thanks[footnoteinfo]{This paper was not presented at any IFAC 
meeting. Corresponding author Mohsen Dlala.}

 \author[Abd]{Abdallah BENABDALLAH}\ead{abdallah.benabdallah@ipeis.rnu.tn},    
\author[Moh,Moh2]{Mohsen DLALA}\ead{mohsen.dlala@ipeis.usf.tn},               

 \address[Abd]{Higher Institute of Computer Science and Multimedia, University of Sfax, Sfax, Tunisia}  
 \address[Moh]{Sfax University, IPEIS, Mathematics-Physics Department, Sfax, Tunisia}             
 \address[Moh2]{Department of Mathematics, College of Sciences, Qassim University, Buraidah, Saudi Arabia}             

\begin{keyword}                           
 Cascaded ODE-heat system; Disturbance to state stabilization; Backstepping stabilization; Observer based control; Input-to-state stability; Output feedback design.

\end{keyword}                             

\begin{abstract}                          

In this paper, we analyze the output stabilization problem for cascaded nonlinear ODE with $1-d$ heat diffusion equation affected by both in-domain and boundary perturbations.  We assume that  the only available part of states is the first components of the ODE-subsystem and one boundary of the heat-subsystem.  The particularity of this system is two folds i) it contains a  nonlinear additive term in the ODE-subsystem, and ii) it is affected by both boundary and in-domain perturbations signals.

For such a system, and unlike the existing works, we succeeded to design an output observer-based feedback that guarantees not only asymptotic stabilization result but also a globally {\it disturbance-to-state stabilization} for our cascaded system. The output feedback is designed using an adequate backstepping transformation recently introduced for coupled ODE-heat equations combined with high-gain observer and high-gain controller.  

\end{abstract}

\end{frontmatter}

\section{Introduction}
Control theory has been revolutionized by the notion of input-to-state stability (ISS) introduced by E.D. Sontag for nonlinear finite-dimensional systems in 1989 in order to study nonlinear ODEs.
 The main idea behind the ISS notion is to unify the  Lyapunov approach and input-output concept. For more details, see the survey paper \cite{Sontag08}.

 Recently, a great deal of work has been done to extend the success of ISS for finite-dimensional systems to infinite-dimensional systems and particularly to systems governed by partial differential equations (PDEs) \cite{Prieur20}, \cite{Jacob22} and \cite{Karafilis19}. For PDEs, a major complication in trying to prove the ISS property is the problem of the existence of solutions for a specific class of input signals. Such a problem is well explained for the semi-linear parabolic boundary control problem in  \cite{Sch20}.

 Let's consider an abstract infinite-dimensional control system
 \begin{equation}
     \label{ISSsys}
     \dot x(t)= \mathcal A x(t)+ \mathcal Bu(t),
 \end{equation}
 where the state $x(t)\in \mathcal H$, the input $u(t)\in\mathcal U$, $\mathcal A$ is an unbounded operator that generates a $C_0$-semi-group  on Hilbert space $\mathcal H$, and operator $\mathcal B :\mathcal U\rightarrow \mathcal H$. It is well known, from the point of view of control theory, that the case where the operator $\mathcal B$ is bounded is similar to that of finite dimensions \cite{Sch20}. The tricky case of ISS arises when $\mathcal B$ is unbounded, in which the presence of \textit{complete mild solutions} becomes non-trivial even for linear systems. In such situation, the ISS property is related to the admissibility notion of control systems \cite{Orliz18}. Formally, we say 
that system (\ref{ISSsys}) is input-to-state stable, if for all signals $u(.)$ (in some Hilbert space $\mathcal{U}$) and initial condition $x_0\in \mathcal H$, there exists a continuous (mild) maximal solution $x(.): [0, +\infty)\rightarrow \mathcal H$  satisfying 
 \begin{equation}\label{ISSest}
     \|x(t)\|_{\mathcal H}\leq \beta(\|x_0\|_{\mathcal H}, t)+\gamma\left( \sup_{s\in[0, t]}\|u(s)\|_{\mathcal U}\right), \; \forall t\geq 0,
 \end{equation}
where $\beta$ is a $\mathcal K\mathcal L$-function, i.e. increasing with respect to the first argument and decreasing to 0 at $+\infty$ with respect to the second argument, and $\gamma$ is a class $\mathcal K$-function, i.e. continuous and increasing with value 0 at zero. Estimation (\ref{ISSest}) can be understood as a stability property of the system (\ref{ISSsys}) with input $(x_0, u)$ and output $x(.)$. Note that in the literature, many definitions and variants of ISS are possible \cite{Prieur20} and in  \cite{Orliz18} for Orlicz spaces which are a natural extension of $L^p$-ISS notion.

Due to the fact that parabolic PDEs are more tractable,  many of the works on ISS property are established for parabolic semi-linear PDEs \cite{Zheng18}, \cite{Mironchenko19}, \cite{Jacoub19} and \cite{Zheng20}. For  hyperbolic equation a sufficient condition for ISS for quasilinear hyperbolic systems with boundary input disturbances is derived in \cite{Bastin21}. Besides, based on Fourier series decomposition, exponential ISS is established for a damped unclamped string with respect to distributed and boundary disturbances is established in \cite{Lhachemi20}.  For coupled ODE-PDE equations, the ISS property has been recently treated in \cite{Lhachemi22} and \cite{ZhangCDC21}  for parabolic coupling, and \cite{Zhang201}
 and \cite{Zhang20} for coupling with wave equation.  Using sliding mode control and the backstepping approach, ISS property is established for coupled parabolic PDEs in \cite{Wang21}. The recently published book \cite{Karafyllis19} provides a detailed study of the ISS property for different classes of PDEs with global Lipchitz nonlinearities and linear differential parts.

Based on our works \cite{Abd20}, \cite{AbdMoh21}, and \cite{Mohsen22} where backstepping state feedback is designed for the coupled nonlinear ODE-heat equation, and inspired by the recent paper \cite{Chang21}, in which an output feedback design is drawn, we deal, in this paper, with observer-based feedback synthesis for cascaded nonlinear ODE with reaction-diffusion equation guaranteeing the ISS against external disturbances. The main novelty compared to previous works is that our observer-based feedback for cascaded nonlinear ODE heat diffusion renders the {\bf closed-loop system  ISS with respect to external perturbations}. In the infinite-dimensional framework, despite constructing a linear output stabilizer for a linear system, the input-to-state stability (ISS) property is guaranteed only for disturbance signals that are admissible with respect to the nominal closed-loop semigroup.

We reach our goal in three main steps. In the first step, we design in section \ref{sec2} a simple Luenberger observer and we use the classical backstepping design for ODE-PDE coupled systems to build our observer-based feedback. In the second step, we prove in section \ref{sec3} the existence of mild solution of the closed-loop system based on Bari Theorem \cite{Urata63}. The third step is devoted to classical  Lypunov analysis which is carried out in section \ref{sec4}. 
In Section  \ref{sec5}, the results of this study are summarized and some future directions are suggested. Finally, some proofs are collected in the Appendix.

\section{Problem statement and controller design}\label{sec2}

We begin our study, in this section, by introducing the cascaded system under consideration and detailing the main assumptions that will be taken into account.

\subsection{Problem formulation}
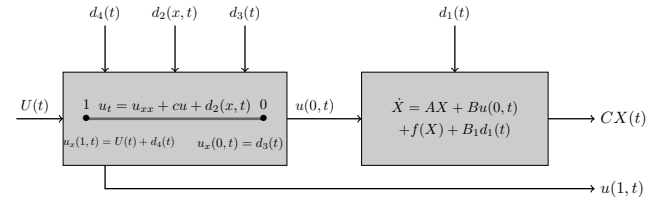
\begin{figure}[!h]
\begin{center}
\scalebox{0.62}{
\begin{tikzpicture}
\draw[thick,->] (-6.4,0) -- (-5.4,0);
\fill[gray!40!white](-5.4,-1) rectangle (-0.6,1);
\draw (-5.4,-1) rectangle (-0.6,1);
\draw[thick,->] (-0.6,0) -- (1,0);
\draw[thick,->] (-4.5,2) -- (-4.5,1);
\draw (-4.5,2) node[scale=0.9, above] {$d_4(t)$};
\draw[thick,->] (-1.5,2) -- (-1.5,1);
\draw (-1.5,2) node[scale=0.9, above] {$d_3(t)$};
\draw[thick,->] (-3,2) -- (-3,1);
\draw (-3,2) node[scale=0.9, above] {$d_2(x, t)$};
\draw[thick,->] (3,2) -- (3,1);
\draw (3,2) node[scale=0.9, above] {$d_1( t)$};
\fill[gray!40!white] (1,-1) rectangle (5,1);
\draw (1,-1) rectangle (5,1);
\draw (-3,0) node[scale=0.9, above]{$u_t=u_{xx}+cu+d_2(x, t)$};
\draw [ultra thick, black!60!white] (-4.9,0) -- (-1.1,0);
\draw (-4.9,0.1) node[scale=0.9, above] {$1$};
\draw (-4.2,-0.3) node[scale=0.7, below] {$u_x(1,t)=U(t)+d_4(t)$};
\draw (-6,0) node[scale=0.9, above] {$U(t)$};
\draw (0,0) node[scale=0.9, above] {$u(0,t)$};
\draw (-4.9,0) node {$\bullet$};
\draw (-1.1,0.1) node[scale=0.9, above] {$0$};
\draw (-1.1,0) node {$\bullet$};
\draw (-1.6,-0.3) node[scale=0.8, below] {$u_x(0,t)=d_3(t)$};
\draw[thick,->] (-4.5,-1) -- (-4.5,-1.5)--(6, -1.5);
\draw (6,-1.5) node[right]{$u(1, t)$};
\draw[thick,->] (5,0) -- (6,0);
\draw (6,0) node[right]{$CX(t)$};
\draw (3,0.25) node[scale=0.9] {$\dot X=AX+Bu(0,t)$};
\draw(3,-0.25) node[scale=0.9]{$+f(X)+B_1d_1(t)$};
\end{tikzpicture}  
}
\end{center}
\caption{The block diagram of the cascaded PDE-ODE system (\ref{sysdimfini})-(\ref{boundlcont}) to boundary control disturbance with output $(CX(t), u(1, t))$. }
\label{fig:PDEODE}
\end{figure}

In this work, we consider a nonlinear ODE cascaded by a $1-d$ reaction-diffusion equation (see Figure \ref{fig:PDEODE}) with a distributed disturbance as well as   disturbances affecting the boundaries.  This system can be represented by the following dynamics:
\begin{eqnarray}
\dot{X}(t)&=& AX(t) +Bu(0,t)+f(X(t))+B_1d_1(t),\label{sysdimfini}\\
u_t(x,t)&=& u_{xx}(x,t)+cu(x, t)+d_2(x, t), \; x\in (0, 1) ,\label{sysheat}\\
u_x(0, t) &=&d_3(t),\label{bound0}\\
u_x(1, t) &=& U(t)+d_4(t),\label{boundlcont}
\end{eqnarray}
where $U(t)\in \mathbb R$ is the control input,  $X(t) \in \mathbb R^n$ is the state of nonlinear ODE subsystem and  $u(x, t) \in \mathbb R$ is the state of the  $1-d$ reaction-diffusion subsystems, respectively.  
The matrices $A \in \mathbb R^{n\times n}$  and $B\in \mathbb R^{n\times 1}$ are given by 
\begin{equation}\label{AB}
A=\left[\begin{array}{ccccc}
0&1&0&\cdots&0\\
0&0&1&\ddots&\vdots\\
\vdots&&\ddots&\ddots&0\\
\vdots&&&\ddots&1\\
0&0&\cdots&\cdots&0\\
\end{array}\right], \;\; B=\left[\begin{array}{c}
0\\
\vdots\\
\vdots\\
0\\
1\\
\end{array}\right], 
\end{equation}
 $c\in \mathbb R$ is scalar constant and $B_1\in \mathbb R^{n\times 1}$ is arbitrarily unknown vector matrices. Distributed disturbances are represented by $d_1(.)$ and $d_2(.)$ in equations (\ref{sysdimfini}) and (\ref{sysheat}), while $d_3(.)$ and $d_4(.)$ are boundary disturbances.

Furthermore, we assume that the only available values of the state $(X(t), u(x, t))$ are 
\begin{equation}\label{output}
y(t)= (CX(t), u(1, t))=(X_1(t), u(1, t)),
\end{equation}
where  $C^T=[1, 0, \cdots, 0], X^T=[X_1, \cdots, X_n]\in \mathbb R^{1\times n}$. For Coupled ODE-PDE systems, it is well known that the stabilization problem becomes very complicated if only a part of state is known. In this setting, besides assumption (\ref{output}) and the presence of non linearity term $f(X)$  in 
equation (\ref{sysdimfini}), we show  that it is still possible to not only obtain asymptotic stability of the cascaded system (\ref{sysdimfini})-(\ref{boundlcont}) but also obtain the ISS stabilization with respect to all disturbances which is called disturbance to state stabilization \cite{Prieur20}.

\subsection{Assumptions and discussion}
Let's introduce the main assumptions taken on   disturbances $d_i(t), i=1,..,4$ and on the nonlinear term $f(X)$ of the ODE subsystem.

\begin{assum} \label{ass1} The disturbances signal $d_i(t) \in L^{\infty}_{\text{loc}}(0, +\infty)$, for $i=1, 3, 4$, i.e. 
\[
\forall i\in\{1, 3, 4\}, \sup_{s\in[0,t]}|d_i(s)| < \infty, \text{ for all } t\geq 0.
\]
\end{assum}

\begin{assum}\label{ass2}The disturbance signal $d_2(., .)$ satisfies the following conditions

\begin{itemize}
    \item For all time $t\geq 0$, $d_2(., t) \in L^2(0, 1)$.
    \item For  almost every $x\in (0, 1)$, $d_2(x, .)\in L^{\infty}_{loc}(0, +\infty).$
    i.e. 
$$\max_{s\in[0,t]}|d_2(x,s)| < \infty,~\text{for all}~t\geq 0.$$
\end{itemize}
\end{assum}
\begin{assum}\label{ass3} There exists a known positive constant $\theta$, such that any $i$-th component of the nonlinear term $f(X)$ satisfies 
\begin{equation}
\left| f_i(X)\right|\leq \theta \sum_{j=i+2}^n|X_j|, \; \; \forall X\in \mathbb R^n, \; 1\leq i\leq n-2,
\end{equation}
and $f_{n-1}(X)=f_n(X)=0$.
\end{assum}

In the following, two remarks should be noted. One is about the type of system being considered, and the second is about perturbations and the nonlinear term.

\begin{rem}
First of all, the class system  (\ref{sysdimfini})-(\ref{boundlcont}) is a well-known and well-studied abstract ODE cascaded with a heat diffusion system.  After being presented as practical models of engineering and physics problems, this type of equation has been formulated and analyzed in \cite{Kris09heat} and \cite{Tang11SCL}. 
Since then, many research papers dealing with the controllability and stabilization of these ODE-heat systems have been published \cite{Abd20}, \cite{Chang21}, \cite{Mohsen22}, \cite{Xie20}, and 
 \cite{Xie21}. Recently, many works analyzing the ISS stability for PDEs systems have been initiated \cite{Prieur20} and \cite{Jacoub19}. 

It is only recently that results have appeared on ISS for coupled parabolic systems \cite{zheng2018iss}, \cite{Karf18}, \cite{Wang21} and \cite{Zhang21}. Regarding this type of coupled systems, the first works on the study of their ISS stability property appeared just recently in \cite{zheng2018iss} and \cite{Karf18}. The ISS of several types of coupled ODE-PDE or PDE-PDE has been  proved in \cite{Karf18}. In \cite{zheng2018iss}, ISS property of coupled beam-string system, was presented.   In \cite{Zhang21}, the ISS property is studied for the linear coupled ODE-heat, while in \cite{Wang21}, a system of coupled linear parabolic PDEs  with memories terms is shown to be ISS stabilizable with respect to disturbances by means of state feedback.

  Recently, it has been established in \cite{dashkovskiy2023robust} an important ISS results for coupled ODE-PDE systems with nonlinearities in both ODE and PDE parts, where disturbances enter through the boundary conditions of the PDE. Their work allows for more general nonlinear structures but addresses only the state feedback case. In contrast, our work focuses on the more challenging output feedback problem, which necessitates additional structural assumptions (specifically, the upper triangular form of nonlinearity in the ODE subsystem) to enable observer design. While seemingly more restrictive, these structural conditions are fundamentally necessary for observer-based control design in nonlinear systems. The distinction between these approaches highlights an important trade-off in nonlinear control theory: more general nonlinear structures are possible when full state information is available, while output feedback requires more specific structural properties to guarantee observability and allow for systematic observer design. An interesting direction for future research would be to investigate whether our explicit output feedback design could be extended to broader classes of nonlinearities while maintaining the observer-based structure, possibly by employing more sophisticated observer designs or alternative backstepping transformations. This could help bridge the gap between the general state feedback results of \cite{dashkovskiy2023robust}  and practical output feedback implementations required in many applications. According to our knowledge, this is the first ISS {\bf{output stabilization}} result  involving a {\bf nonlinear} term in the ODE subsystem for this class of coupled ODE-heat systems. 
\end{rem}

\vspace{3mm}

\begin{rem}
Under Assumption \ref{ass3}, the nonlinear term $f(X)$ in the ODE-subsystem (\ref{sysdimfini}) has an upper triangular structure called feedforward nonlinear systems  \cite{Mazenc1996adding} which are well studied in the literature. 
This triangular structure allows us to design an observer-based feedback and to avoid the strong global Lipschitz assumption on the nonlinear term as assumed for example in \cite{Ahmed15}.  To succeed in this task, our method relies on two paths, we adopt the classical technique of combining a high-gain observer with a high-gain controller, as shown for example in \cite{Hao07}, for the ODE part and then we use backstepping design method \cite{Tang11} and \cite{Kris09heat} for coupled ODE-heat equation to elaborate the output feedback for the whole cascaded system.
\end{rem}

In the following, we write system (\ref{sysdimfini})-(\ref{boundlcont}) in the abstract semigroup form and then we establish the existence and uniqueness of maximal mild solutions.

\section{Well-posedness }\label{sec3}
First, we introduce the state space of the solutions for the cascaded system (\ref{sysdimfini})-(\ref{boundlcont}). We denote $L^2(0, 1)$ the Hilbert space consisting of measurable square  integrable functions on $(0, 1)$ being equipped by its canonical inner product $<\cdot, \cdot>_{L^2(0, 1)}$ and its induced norm $\|\cdot\|$. Moreover, we consider the energy Hilbert space $\mathcal H=\mathbb R^n\times L^2(0, 1)$ and its canonical inner product 
$$\left< (X, u), (Y, v)\right>_{\mathcal H} = X^TY+<u, v>_{L^2(0, 1)},$$
for all $(X, u), (Y, v)\in\mathcal H$. In addition, consider the Hilbert space $\mathcal H_1=\mathbb R^n\times \mathcal D$, where 
$$\mathcal D=\left\{u \in H^2(0, 1), \; u'(0)= u'(1)=0\right\}.$$
 Define operator $\mathcal A$ as follows:
 \begin{eqnarray}
 \mathcal A(X, u)&=&(AX +B u(0), u''+c u),\; \forall (X, u)\in \mathcal D(\mathcal A),
 \end{eqnarray}
 where $\mathcal D(\mathcal A)=\mathcal H_1$.
 Then, system (\ref{sysdimfini})-(\ref{boundlcont}) can be written in the following compact form
 \begin{eqnarray} 
 \label{dotX=AX++}
 \dot{\mathcal X}(t)&=&\mathcal A \mathcal X(t)+\mathcal F(\mathcal X(t))+\mathcal  B_1 U(t)+ \mathcal B_2 \tilde d_1(t)\nonumber\\
 &&+\mathcal B_3 d_3(t)+\mathcal B_4 d_4(t),\\
 \mathcal Y(t)&=&\mathcal C \mathcal X(t),\label{outputYt}
 \end{eqnarray}
where $\mathcal X(t)=\left[X(t)^T, u(\cdot, t)\right]^T$, 
$\mathcal F(\mathcal X(t))=[f(X(t))^T, 0]^T$, $\mathcal B_1=[0_{\mathbb R^n}^T, \delta(x-1)]^T$, $\mathcal B_2=\left[\begin{array}{cc} B_1& 0_{\mathbb R^n}\\
0& 1
\end{array}\right]$, $\mathcal B_3=\left[ 0_{\mathbb R^n}^T, \delta(x)\right]^T$, $\mathcal B_4=\left[ 0_{\mathbb R^n}^T, \delta (x-1)\right]^T, \tilde d_1(t)=[d_1(t), d_2(\cdot, t)]$ and $\mathcal C=[C, \delta(x)]$. Here $\delta(.)$ denotes the Dirac function.

We point out, that establishing disturbance-to-state stabilization for the system (\ref{dotX=AX++}) is not trivial since disturbance operators $\mathcal B_1,$ $\mathcal B_3,$ and $\mathcal B_4$ are unbounded.
Now, we state the first main result concerning the well-posedness of system (\ref{dotX=AX++}). 

\begin{thm}\label{Th1} Under Assumptions \ref{ass1}, \ref{ass2} and \ref{ass3}, system (\ref{dotX=AX++}) is well-posed, i.e. for any initial condition $\mathcal X_0=(X_0, u_0)\in \mathcal H$, and $U(\cdot), \tilde d_1, d_3(\cdot)$ and $d_4(\cdot) \in L^2_{loc}(0, +\infty)$, system (\ref{dotX=AX++}) admits a unique mild solution 
\begin{align}
\mathcal X(t)=\left(X(t), u(\cdot, t)\right) \in C(0, T_{\max}, \mathcal H).
\end{align} 
 Moreover,  the following prior estimation is satisfied
\begin{eqnarray}
\|\mathcal X(t)\|_{\mathcal H}&&\leq c_t \|\mathcal X(0)\|_{\mathcal H}+d_t  \int_0^t U^2(s) ds\nonumber\\
&&+d_t\sup_{s\in [0, t]}\left(\|d_2(., s)\|+\sum_{i=1, 3, 4}|d_i(s)|\right),\label{priorest}
\end{eqnarray}
for all $t\in [0, T_{\max})$, where $c_t$ and $d_t$ are time-depending  positive reals.
\end{thm}

The proof of Theorem \ref{Th1}, presented in Appendix \ref{pfTh1}, employs semigroup theory and the Bari Theorem in a novel manner. The application of the Bari Theorem requires careful consideration, as we introduce a modified basis to establish that the family of eigenvectors constitutes a Riesz basis. Furthermore, we demonstrate the admissibility of operators $\mathcal B_1$, $\mathcal B_2$, and $\mathcal B_3$ with respect to the semigroup $T(t)$, following the framework developed in \cite{Sch20}.

\begin{rem}
Prior estimation (\ref{priorest}) shows that the energy of system (\ref{sysdimfini})-(\ref{boundlcont}) is bounded when input signals $(U(.), \tilde d_1(.), d_3(.), d_4(.))$ are bounded.
\end{rem}

Our main objective is to design a dynamic output feedback that rends the system ISS with respect to perturbation $(\tilde d_1(.), d_3(.), d_4(.))$. As we go on, we will look at our second main finding.

\section{Global exponential disturbance-to-state stabilization}\label{sec4}

Under Assumptions 1, 2 and 3, we present a controller design procedure that makes a cascaded system (\ref{sysdimfini})-(\ref{boundlcont}) globally exponential disturbance-to-state stable. Consider the set $\mathcal K$ of all continuous increasing functions $k(.): [0, +\infty)\rightarrow [0, +\infty)$ such that $k(0)=0$. Now, let's introduce the following definition of global exponential disturbance-to-state stabilization by dynamic output feedback given in \cite{Prieur20}.

\begin{defn} System (\ref{sysdimfini})-(\ref{boundlcont}) is said to be globally exponentially disturbance-to-state  stabilizable (DISS) by dynamic output feedback, if the exists a dynamic system
\begin{equation}\label{dynacontrol}
\begin{array}{lll}
\dot \xi(t)&=&h_c(\xi(t), y(t)),\\
U(t)&=&h_i(\xi(t), y(t))
\end{array}
\end{equation}
such that the closed loop system (\ref{sysdimfini})-(\ref{boundlcont}) with (\ref{dynacontrol}) admits  a continuous mild solutions for all initial conditions $(X_0, u_0, \xi_0)$ in some Hilbert space, and there exist   a positive constant $ \omega$ and a class $\mathcal K$ function $\Gamma$ such that
\begin{equation}\label{estimationISS}
\left|X(t)\right|+ \|u(., t)\|\leq cM_{(X_0, u_0, \xi_0)}e^{-\omega t}+\Gamma\left( \sup_{s\in [0, t]}|D(s)|\right),
\end{equation}
where $M_{(X_0, u_0, \xi_0)}$ depending on the initial condition $(X_0, u_0, \xi_0)$ and $$|D(t)|=|d_1(t)|+\|d_2(\cdot, t)\|+|d_3(t)|+|d_4(t)|.$$
\end{defn}

To achieve the desired control objective, we use the indirect method  (terminology introduced in \cite{Praly09}) for output feedback stabilization. In an observer-based approach, this indirect method consists in analyzing the stability properties of the interconnecting system formed by error dynamics and controller dynamics.
This method is recently used for  global exponential stabilization by output feedback of the unperturbed cascaded ODE-PDE system \cite{Chang21}. Indeed, it allows to bypass all the difficulties of the system related to the nonlinear term encountered in \cite{Abd20} after having applied the backstepping transformations. 

Our stabilization strategy is based on the following three successive steps.  First, for the cascaded system (\ref{sysdimfini})-(\ref{boundlcont}), we construct a Luenberger observer system $(\mathcal L)$ comprising a copy of the system (\ref{sysdimfini})-(\ref{boundlcont}) and an added corrective term. Second, we build an output feedback $U(t)$ for system $(\mathcal L)$  using backstepping method for coupled ODE-PDE. Finally, based on the Lyapunov method, we demonstrate that the coupled system consists of the error dynamics and the observer system $(\mathcal L)$ in closed loop with an output feedback $U(t)$ is  globally exponentially ISS with respect disturbances signals.

\subsection{Observer design}
In this section, we turn to estimate the unavailable system state $(X(t), u(x,t))$ based on the measured output $y(t)$. In order to achieve this, the observer $(\mathcal L)$ is traditionally introduced in the following way:
\begin{eqnarray}
\dot{\hat{X}}(t)&=& A\hat X(t) +B\hat u(0,t)-D_rLC\tilde X(t),\label{obsdimfini}\\
\hat u_t&=& \hat u_{xx}+c\hat u+k(x)\tilde u(1, t), \; x\in (0, 1) ,\label{obsheat}\\
\hat u_x(0, t) &=&0,\label{obsbound0}\\
\hat u_x(1, t) &=& U(t)+q_2\tilde u(1, t),\label{obsboundlcont}
\end{eqnarray}
where $L\in\mathbb R^{n\times 1}, k(\cdot), q_2, r\geq 1$ are observer gains to be designed, the diagonal matrix $D_r=diag(r^{-1}, \cdots, r^{-n})$,  $\tilde u(1, t)=u(1, t)-\hat u(1, t)$ and $\tilde X(t)=X(t)-\hat X(t)$. It is important to emphasize that system (\ref{obsdimfini})-(\ref{obsboundlcont}) corresponds to the Luenberger observer, often referred to as a replica of the nominal system, but with an additional correction term included. Evidently, system (\ref{obsdimfini})-(\ref{obsboundlcont}) is both free from disturbances and linear. As a result, constructing a stabilizing feedback for this system becomes more straightforward, and we will subsequently demonstrate that this stabilizing feedback is also effective in stabilizing the system (\ref{sysdimfini})-(\ref{boundlcont}). This control design strategy is widely used in the control theory community.  In the following, we use  backstepping method for coupled ODE-PDE systems. To be more precise, we will build
an output feedback that guarantees a global exponential disturbance-to-state stabilization for couples systems (\ref{sysdimfini})-(\ref{boundlcont}) via observer ($\mathcal L$). 

\subsection{Backstepping controller design}

In this section, the backstepping transformation is implemented  based on the work of \cite{Kris09heat} and \cite{Tang11} (See \cite{AbdMoh21} and \cite{Mohsen22} for a closer transformation). Let $r$ be a positive constant gain and $\Delta_r = diag( r^{-n},   r^{-n+1},\cdots, r^{-1})$ a diagonal matrix. Consider the following  backstepping transformation  
\begin{eqnarray}
\hat Z(t)&=& \Delta_r \hat X(t),\label{transXZ}\\
\hat w(x, t)&=& \hat u(x,\!t)\!-\!\int_0^xs(x,\!y)\hat u(y,\! t)dy-\psi(x)\Delta_r\hat{X}(t),\label{transwu}
\end{eqnarray}
where  the  kernel $s(x, y)$ 
satisfies the following equations
\begin{eqnarray}
s_{xx}(x, y)-s_{yy}(x, y)&=& c s(x, y),\; \forall (x, y) \in \mathcal T_1,\label{sxx1}\\
s(x, x)&=&-\frac{c}{2}x, \; \forall x\in [0, 1],\label{sxx2}\\
s_y(x, 0)&=&-r^{-1}\psi(x)B,  \; \forall x\in [0, 1].\label{sxx3}
\end{eqnarray}
And the function $\psi(x)$ satisfies
\begin{eqnarray}
\psi''(x)-r^{-1}\psi(x)A&=& 0,\; \forall x, y \in \mathcal T_1,\label{psi1}\\
\psi(0)&=&K, \label{psi2}\\
\psi'(0)&=&0.\label{psi3}
\end{eqnarray}
 The following lemma proves the existence and gives the explicit expressions of functions $s(x, y)$  and $\psi(x)$.

\begin{lem}

\begin{enumerate}
  \item   For all $x\in[0,1]$, we have
\begin{align}\label{psi}
    \psi(x)=K\sum_{k=0}^{\infty} \frac{x^{2k}}{(2k)!}(r^{-1}A)^k
\end{align}
Moreover, there exists a positive constant $c_{\psi}$ such that,  for all $x\in [0, 1]$,
\begin{equation}\label{cpsi}
    |\psi(x)|\leq c_{\psi}.
\end{equation}
 \item  System (\ref{sxx1})-(\ref{sxx3})  admits  a unique solution $s(x, y)$ in $C^2(\mathcal T_1)$. In addition, the transformation (\ref{transXZ})-(\ref{transwu}) is invertible and its inverse is linearly bounded.
\end{enumerate}
\end{lem}

\begin{pf}
    \begin{enumerate}
        \item The proof of the first item can be found in \cite{Lei20}.
\item The proof of the second item is a particular case of the proof of Proposition 2.1 in \cite{Zhang21}.
\end{enumerate}
    \end{pf}
By applying transformation (\ref{transXZ}) and (\ref{transwu}),  system (\ref{obsdimfini})-(\ref{obsboundlcont}) is converted into 
\begin{align}
r\dot{\hat{Z}}(t)\!&=\!(A\!+\!BK)\hat Z(t)\! +\!B\hat w(0,\!t)\!-\!r^{-n}\!L\! C\! \tilde X^1(t),\label{bacobsdimfini}\\
\hat w_t(x,t)&=\!\hat w_{xx}(x,t)\!+\!\tilde k(x)\tilde u(1, t)\!+\!\frac{\psi(x)}{r^{n+1}}L C\tilde X^1(t),\label{bacobsheat}\\
\hat w_x(0, t) &=0,\label{bacobsbound0}\\
\hat w_x(1, t) &=-\eta \hat w(1, t),\label{scaledobsboundlcont}
\end{align}
where  the output feedback $U(t)$ is chosen as
\begin{eqnarray}
U(t)&=&\left(q_2-\eta-c/2\right)\hat u(1, t)-q_2u(1, t)\notag\\&&+\int_0^1\left(s_x(1, y)+\eta s(1, y)\right)\hat u(y, t)dy\nonumber\\
&&+(\psi'(1)+\eta\psi(1))\Delta_r \hat X(t).\label{outputfeedback}
\end{eqnarray}

and the function $\tilde k(x)$ is given by
$$\tilde k(x)=k(x)-\int_0^xs(x, y)k(y)dy, \; \forall x\in [0, 1].$$

Now, we demonstrate that by suitably manipulating the gain matrices $K$ and $L$, the kernel $k(\cdot)$, and the gain scalars $r$, $\eta$, and $q_2$, the output feedback (\ref{outputfeedback}), together with (\ref{obsdimfini})-(\ref{obsboundlcont}), enables the system (\ref{sysdimfini})-(\ref{boundlcont}) to achieve global disturbance-to-state stability with respect to disturbances. To accomplish this, we must carefully choose an appropriate Lyapunov function for the closed-loop system in the subsequent subsection.

\subsection{Target closed-loop system}

First, we introduce the two components of the scaled state error as follows:
\begin{align}
\tilde X(t) &= \Delta_r\left(X(t)-\hat X(t)\right) \label{errortrans1} \\
\tilde u(x, t) &= u(x, t) - \hat u(x, t). \label{errortrans2}
\end{align}

Next, upon performing a straightforward computation, we find that the dynamics governing $(\tilde X(t), \tilde u(x, t))$ are given by
\begin{align}
r\dot{\tilde{X}}(t)&= (A+LC)\tilde X(t) +B\tilde u(0,t)+r\Delta_rf(X(t))\notag\\
&+r\Delta_rB_1d_1(t),\label{errordimfinitargdisdis}\\
\tilde u_t(x,t)&\!=\! \tilde u_{xx}(x,\!t)\!+\!c\tilde u(x,\!t)\!-\!k(x)\tilde u(1,\!t)\!+\!d_2(x,\!t),\label{errorheattargdisdis}\\
\tilde u_x(0, t) &=d_3(t),\label{errorbound0targdisdis}\\
\tilde u_x(1, t) &=-q_2\tilde u(1, t)+d_4(t).\label{errorboundlconttargdisdis}
\end{align} 
To derive a stabilizing term for the nominal part of equation (\ref{errorheattargdisdis}) (similar to the approach used for subsystem (\ref{errordimfinitargdisdis})), we initially employ the following classical transformation
\begin{equation}\label{transuw}
\tilde u(x, t)\!=\!\tilde w(x,\!t)\!-\;\int_x^1q(x,\!y) \tilde w(y, t)dy, \; \forall x\in [0, 1],
\end{equation}
where the kernel $q(x, y)$ is selected as usual  satisfying (see \cite{Kris04})
\begin{eqnarray}
q_{xx}(x, y)-q_{yy}(x, y)&=&-c q(x, y),\label{qxx1}\\
q(x, x)&=&-\frac{cx}{2},\label{qxx2}\\
q_x(0, y)&=&0,\label{qxx3}
\end{eqnarray}
for    all $(x, y)\in \mathcal T_2=\{(x, y)\in [0, 1]^2, \; x\leq y\leq 1\}$.
The explicit expression of kernel $q(x, y)$  solution of (\ref{qxx1})-(\ref{qxx3})  is given  by 
(see \cite{Kris04})
\begin{equation}\label{q}
q(x, y)=-cy\frac{I_1\left(\sqrt{c(y^2-x^2)}\right)}{\sqrt{c(y^2-x^2)}}, \; \forall (x, y)\in \mathcal T_2,
\end{equation}
where $I_1$ is the first kind of modified Bessel function. 

Existence and properties of the inverse transformation of (\ref{transuw}) are given in the following lemma.

\begin{lem} \label{invtransuw}Transformation (\ref{transuw}) is invertible and 
\begin{equation}\label{transuwinv}
\tilde w(x,\!t)=\tilde u(x,\!t)\!-\!\int_x^1p(x,\!y) \tilde u(y,\!t)dy, \; \forall x\in [0, 1],
\end{equation}
where the kernel $p(x, y)$ is given by
\begin{eqnarray}
p_{xx}(x, y)-p_{yy}(x, y)&=&-c p(x, y),\label{pxx1}\\
p(x, x)&=&-\frac{cx}{2},\label{pxx2}\\
p_x(0, y)&=&0,\label{pxx3}
\end{eqnarray}
for all $(x, y) \in \mathcal T_2$.
\end{lem}
The proof of the Lemma \ref{invtransuw} is given in Appendix \ref{proofinvtransuw}.

Transformation (\ref{transuw}) allows us to obtain the cascade target system composed of both stable ODEs and PDEs as follows:
\begin{eqnarray}
r\dot{\tilde{X}}(t)&=&\!(A\!+\!LC)\tilde X(t)\!+\!B\tilde w(0,t)\!-\!B\int_0^1\!q(0,\!y)\tilde w(y,\!t)dy\nonumber\\&&+r\Delta_rf(X(t))+r\Delta_r B_1d_1(t),\label{errordimfinitargtarg}\\
\tilde w_t(x,t)&=& \tilde w_{xx}(x,t)+\tilde k(x)\tilde w(1, t)+\tilde d_2(x, t),\label{errorheattargtarg}\\
\tilde w_x(0, t) &=&d_3(t),\label{errorbound0targtarg}\\
\tilde w_x(1, t) &=&-(q_2-c/2)\tilde w(1, t) +d_4(t),\label{errorboundlconttargtarg}
\end{eqnarray}
where  $\tilde d_2(x, t)$  and $\tilde k(x)$ are given by
\begin{align}
    \tilde k(x)&=-k(x)+p(1, 1) q_2+p_y(x, 1)\notag\\
    &-\int_x^1p(x, y) k(y)dy, \label{tildek}\\
\tilde d_2(x, t)&= d_2(x, t)-p(1, 1)d_4(t)\notag\\
&-\int_x^1p(x, y) d_2(y, t)dy,\label{tilded2}
\end{align}
for all $x\in [0, 1]$ and $t\geq 0$.

At this stage, we can design the observer gain $k(x)$ in the following way.

\begin{lem}\label{lemk} Let $h : [0, 1]\rightarrow \mathbb R$ be a continuous function. Then, there exists a gain $k(x)$ solution of the following integral equation 
\begin{equation}\label{k}
    k(x)=h(x)-\int_x^1p(x, y) k(y) dy, \;\forall x\in [0, 1].
\end{equation}
\end{lem}
The proof is based on the standard successive approximation and can be found in Appendix \ref{pflemk}. 

Based on Lemma \ref{lemk}, we select the gain $k(x)$ as 
\begin{align}\label{gaink}
    k(x)=p(1, 1) q_2+p_y(x, 1)-\int_x^1p(x, y) k(y)dy,
\end{align}
 in such a way that $\tilde k(x)$ in (\ref{tildek}) is vanishing, i.e. $\tilde k(x)=0$ for all $x\in [0, 1]$. 

Once all the above transformations have been carried out and  the kernel $k(x)$ is chosen as in (\ref{gaink}), the foundations have been laid for the Lyapunov analysis of the closed-loop system, which will be the subject of the next section.

\subsection{Stability analysis of the closed-loop system}

In light of the fact that $\tilde k=0$, let's write the  closed-loop system (\ref{errordimfinitargtarg})-(\ref{errorboundlconttargtarg}) with (\ref{bacobsdimfini})-(\ref{scaledobsboundlcont}) as follows:
\begin{eqnarray}
r\dot{\hat{Z}}(t)&=& \left(A+BK\right)\hat Z(t) +B\hat w(0,t)-LC\tilde X(t),\label{closeddimfin1}\\
 r\dot{\tilde{X}}(t)&=& (A+LC)\tilde X(t) +B\tilde w(0,t)\nonumber\\&&-B\int_0^1q(0, y)\tilde w(y, t)dy+r^{-1}\Delta_rf(X(t))\nonumber\\
 &&+r^{-1}\Delta_r B_1d_1(t),\label{closeddimfin2}\\
 \hat w_t(x,t)&=& \hat w_{xx}(x,t)+r^{-1}\psi(x)LC\tilde X(t), \; x\in (0, 1),\label{closedheat1}\\
\tilde w_t(x,t)&=& \tilde w_{xx}(x,t)+\tilde d_2(x, t), \label{closedheat2}\\
\hat w_x(0, t) &=&0,\label{closedbound1}\\
\tilde w_x(0, t)&=&d_3(t),\label{closedbound2}\\
\hat w_x(1, t) &= &-\eta \hat w(1, t),\label{closedbound11}\\
\tilde w_x(1, t) &=&-(q_2-c/2)\tilde w(1, t) +d_4(t),\label{closedbound22}
\end{eqnarray}
for all $x\in (0, 1)$ and $t\in [0, T_{\max}).$  Our second main result is summarized in the following Theorem which states that the closed-loop system (\ref{closeddimfin1})-(\ref{closedbound22}) satisfies DISS estimation (\ref{estimationISS}).
\begin{thm}\label{isstheorem} 
 Let assumptions \ref{ass1}, \ref{ass2} and \ref{ass3} from Theorem \ref{Th1} be satisfied. Then,

 \begin{enumerate}
     \item[(i)] For any initial conditions $(X(0),u(.,0)) \in \mathcal D(\mathcal A)$ and observer initial conditions $(\hat X(0), \hat u(.,0)) \in \mathcal H$, the closed-loop system (\ref{closeddimfin1})-(\ref{closedbound22}) admits a unique solution with

$$(X, u, \hat X, \hat u) \in  C([0, \infty), \mathcal H \times \mathcal H).$$ 
\item[(ii)]  The system is ISS in the following sense: there exist $\beta\in  \mathcal K\mathcal L$ and $\gamma\in \mathcal K$ such that for all $t \geq 0$:
 \begin{align}\label{ISSthm}
 & \| (X(t), u(.,t)) \|_{\mathcal H}+\| (\hat X(t), \hat u(.,t)) \|_{\mathcal H}\leq  \notag\\ &\beta\left(\|(X(0),u(.,0), \hat X(0), \hat u(.,0)) \|_{\mathcal H\times \mathcal H}, t\right)+\gamma\left(\sup_{s\in [0, t]}D(s)\right),\notag\\&
  \end{align}
 for all $t\geq 0$, where $D(.) $ represents the collection of all disturbance functions defined in  (\ref{Dt}).
 \end{enumerate}
\end{thm}

\begin{pf}\label{isstheoremprf} 
Three main steps make up the proof.
We derive two Lyapunov functions in the first two steps, a first one for subsystems 
(\ref{closeddimfin1}), (\ref{closedheat1}), (\ref{closedbound1}), and (\ref{closedbound11}), and a second one  is established for subsystems (\ref{closeddimfin2}), (\ref{closedheat2}), (\ref{closedbound2}), and (\ref{closedbound22}).
In the third step, utilizing a composite Lyapunov function, we establish the global exponential input-to-state stability of the whole system (\ref{closeddimfin1})-(\ref{closedbound22}) with respect to all disturbance signals.

\noindent\textbf{First step :} \textbf{Lyapunov function for subsystem} (\ref{closeddimfin1}), (\ref{closedheat1}), (\ref{closedbound1}), (\ref{closedbound11}). 

Utilizing the standard pole placement theorem for linear control systems and selecting any positive constant $\delta_1$, it is possible to find a gain matrix $K\in\mathbb R^{1\times n}$ and a positive definite matrix $P_1=P_1^T>0$ that fulfill the following matrix inequality: 
\begin{equation}\label{lyapP1}
P_1(A+BK)+(A+BK)^TP_1\leq -\delta_1 I_n.
\end{equation}
It is noteworthy that through a straightforward manipulation, the bilinear matrix inequality (\ref{lyapP1}) with unknowns $P_1$, $K$, and $\delta_1$ can be reformulated as a linear matrix inequality, enabling its numerical solution (see \cite{Boyd94}). This computational transformation is essential in explicitly computing the kernel $\psi(\cdot)$ of the output feedback (\ref{outputfeedback}).
 Let $a$ be a positive constant to be selected adequately, and consider the following Lyapunov function $V_1$  for the scaled observer system (\ref{closeddimfin1}), (\ref{closedheat1}), (\ref{closedbound1}), (\ref{closedbound11}) 
\begin{equation}\label{V11}
    V_1(t)=\hat Z(t)^TP_1 \hat Z(t)+\frac{a}{2}\int_0^1\hat w(x, t)^2dx.
    \end{equation}
        
 The properties of the Lyapunov function $V_1(t)$ are given by the following lemma. 

 \begin{lem} \label{V1} Let $a>2$.
  Then, the Lyapunov function $V_1(t)$ defined in (\ref{V11}) satisfies the following properties:
\begin{enumerate}
    \item There exist two positive constants $\nu_1$ and $\nu_2$ such that
    \begin{equation}\label{V1q}
        \nu_1\left\|(\hat Z, \hat w)\right\|_{\mathcal H}^2\leq V_1\leq \nu_2\left\|(\hat Z, \hat w)\right\|_{\mathcal H}^2.
        \end{equation}
     \item The time derivative of $V_1(t)$ along the trajectories of (\ref{closeddimfin1}), (\ref{closedheat1}), (\ref{closedbound1}), (\ref{closedbound11}) satisfies
    \begin{align}
\dot V_1(t)&\leq-r^{-1}\left(\delta_1-r^{-1}|P_1|^2-1\right)|\hat Z(t)|^2\notag\\
&-\left(\frac{a}{4}-\frac{1}{2}-\frac{a}{2}r^{-1}c_{\psi}|L|\right)\|\hat w(., t)\|^2\nonumber\\
& -\left(a(\eta-\frac{1}{2})-1\right)\hat w(1, t)^2\notag\\
&+ \left(r^{-1}|P_1L|^2+\frac{1}{2}ar^{-1}c_{\psi} |L|\right) |\tilde X(t)|^2,\label{dotV112} 
 \end{align}
    for all $t\in [0, T_{\max})$.
\end{enumerate}
\end{lem}   
    
The proof of  Lemma \ref{V1} is so classical and can be found in Appendix \ref{lemm3pf}.

\noindent\textbf{Second step :} \textbf{Lyapunov function for subsystem}  (\ref{closeddimfin2}), (\ref{closedheat2}), (\ref{closedbound2}), (\ref{closedbound22}). 

Similarly, by applying the pole placement theorem \cite{Wonham67}, we can find a gain matrix $L\in\mathbb R^{1\times n}$ and a positive definite matrix $P_2=P_2^T>0$ such that  the following matrix inequality holds true
\begin{equation}\label{lyapP2}
P_2(A+LC)+(A+LC)^TP_2\leq -\delta_2 I_n.
\end{equation}
For the scaled observer (\ref{closeddimfin2}), (\ref{closedheat2}), Let consider the following Lyapunov function 
\begin{equation}\label{V22}
    V_2(t)=\tilde X(t)^TP_2\tilde X(t)+\frac{b}{2}\int_0^1\tilde w(x, t)^2dx,
    \end{equation}
where $b$ is a positive constant to be selected in the sequel. The Lyapunov function $V_2(t)$ satisfies the following estimations. 

\begin{lem} \label{V2} The Lyapunov function $V_2(t)$ defined in (\ref{V22}) satisfies the following.
\begin{enumerate}
    \item There exist two positive constants $\mu_1$ and $\mu_2$ such that
\begin{align}\label{V21}
\mu_1\left(|\tilde X|^2\!+\!\|\tilde w\|^2\right)\leq V_2\leq \mu_2\left(|\tilde X|^2\!+\!\|\tilde w\|^2\right).
\end{align}
    \item The time derivative of $V_2(t)$ along the trajectories of (\ref{closeddimfin2}), (\ref{closedheat2}), (\ref{closedbound2}), (\ref{closedbound22}) satisfies
\begin{align}\label{dotV2121}
\dot V_2(t)&\leq  
-r^{-1}\left(\delta_2- r^{-1} m_2\right) |\tilde X(t)|^2\notag\\
&+(m_1+\frac{5-b}{4})\|\tilde w(., t)\|^2\nonumber\\
&-\left(b(q_2-\frac{c}{2}-\frac{1}{2})-\frac{1}{2}b^2-\frac{3}{2}\right)\tilde w(1, t)^2\notag\\
&+ n\theta\lambda_{\max}(P_2)r^{-2}|\hat Z(t)|^2+d_1^2(t)\nonumber\\
& +\frac{1}{2}b^2d_3(t)^2+\frac{1}{2}d_4(t)^2+\frac{1}{2}b^2\|\tilde d(., t)\|^2,
\end{align}   
where $m_1=\int_0^1q(0, y)^2dy$ and $m_2=\lambda_{\max}(P_2)^2(2+|B_1|^2)+3n \theta\lambda_{\max}(P_2).$
\end{enumerate}
\end{lem}
The proof of Lemma \ref{V2} is omitted to the Appendix \ref{lemma5prf}.


\noindent{\bf Third Step : Composite Lyapunov function. }

To achieve the proof of our Theorem, let's consider the following composite  Lyapunov function
\begin{equation}\label{Lyaptotal}
\mathcal V(t)=V_1(t)+ \gamma V_2(t),
\end{equation}
where $\gamma$ is a positive constant to be selected.
 By (\ref{dotV112}) and (\ref{dotV2121}), we obtain
\begin{eqnarray}\label{dotVV}
\dot{\mathcal V}(t)&\leq& -\tau_1 |\hat Z(t)|^2-\tau_2\|\hat w(., t)\|^2-\tau_3|\tilde X(t)|^2\notag\\
&-&\tau_4\|\tilde w(., t)\|^2-\tau_5\tilde w(1, t)^2-\tau_6\hat w(1, t)^2\notag\\
&+&\gamma d_1^2(t)\!+\!\frac{\gamma}{2}b^2d_3(t)^2\!+\!\frac{\gamma}{2}d_4(t)^2\!+\!\frac{\gamma}{2}b^2\|\tilde d(., t)\|^2,
\end{eqnarray}
where
\begin{eqnarray}
\tau_1&=&r^{-1}\left(\delta_1-r^{-1}|P_1|^2-1-\gamma r^{-1}n\theta\lambda_{\max}(P_2)\right),\\
\tau_2&=&\frac{a}{4}-\frac{1}{2}-\frac{a}{2}r^{-1}c_{\psi}|L|,\\
\tau_3&=&r^{-1}\left(\gamma(\delta_2-|P_1L|^2-r^{-1} m_2) -\frac{1}{2}ac_{\psi} |L|\right),\\
\tau_4&=&\frac{\gamma}{4}\left(b-4m_1-5\right),\\
\tau_5&=&-\frac{\gamma}{2}b^2+b\left(q_2-\frac{c}{2}-\frac{1}{2}\right)-\frac{3}{2},\\
\tau_6&=&a\left(\eta-\frac{1}{2}\right)-1.
\end{eqnarray}
 Now, let $\delta_1>1$ and $r>2m_2/\delta_2$. It follows $\tau_3\geq r^{-1}\left(\gamma\delta_2/2-|P_1L|^2 -\frac{1}{2}ac_{\psi} |L|\right)$. Taking $$\gamma>\gamma^*:=\frac{|P_1L|^2 +ac_{\psi} |L|}{2\delta_2}$$ leads to $\tau_3>0.$

By selecting $r>r^*:=\max\left\{\frac{2 m_2}{\delta_2}, \frac{|P_1|^2+n\gamma\theta\lambda_{\max}(P_2)}{\delta_1-1}, \frac{2 a c_{\psi}|L|}{a-2}\right\}$, the positivity of $\tau_1$ and   $\tau_2$ is guaranteed.

In addition, for $b> 4m_1+5$, we get $\tau_4>0$. On the other hand, the second-degree polynomial in $b$
\begin{align}
g(b)=b^2-\left(2q_2-c-1\right)b+3
\end{align}
 has two positive  roots $b_1<b_2$. Let $\displaystyle b=\frac{b_1+b_2}{2}=q_2-c/2-1/2>4m_1+5$. Thus, to ensure the positivity of $\tau_5$, we select $q_2>q_2^*:=4m_1+\frac{c}{2}+\frac{11}{2}$.
  Finally, by taking $\eta>1$ we can guarantees that
 $\tau_6>0$, and then it yields  that 
 all constants $\tau_i, i=1, \cdots, 6$ are positive.
 
 From (\ref{Lyaptotal}) and by using (\ref{V1q}) and (\ref{V21}) the exists a positive constant $\tau$ such that
 \begin{align}\label{dotmathcal{V}}
 \dot{\mathcal V}(t)\leq -\tau \mathcal V(t) +D(t), \forall t\in [0,T_{\max}),
 \end{align}
 where  
 \begin{align} \label{Dt}
 D(t)=d_1^2(t)\!+\!\frac{1}{2}b^2d_3(t)^2\!+\!\frac{1}{2}d_4(t)^2\!+\!\frac{1}{2}b^2\|\tilde d(., t)\|^2
 \end{align}
   
  Integrating (\ref{dotmathcal{V}}) on $[0,T_{\max})$ gives
  \begin{align}\label{v-t}
 \mathcal  V(t)\leq e^{-\tau t}\mathcal V(0)+\frac{D_{T_{\max}}}{\tau},
  \end{align}  
  where  
 \begin{align} \label{Dmax}
 D_{T_{\max}}=\sup_{t\in[0,T_{\max})}D(t).
 \end{align}
 By the way of contradiction, suppose that $T_{\max}< \infty$.
 Upon examining (\ref{v-t}), it becomes evident that the solution $\mathfrak{X}(t) = \left(\hat{Z}(t), \hat{w}(.,t), \tilde{X}(t), \tilde{w}(.,t)\right)$ remains bounded on the interval $[0, T_{\max})$. However, this contradicts Theorem 1.4  in \cite{Pazy}. Therefore, we can conclude that $T_{\max} = \infty$.
 Now, a simple integration of (\ref{dotmathcal{V}}) gives
  \begin{align}\label{issclosed}
   \mathcal  V(t)\leq e^{-\tau t}\mathcal V(0)+\frac{1}{\tau}\sup_{s\in[0,t]}D(s),
  \end{align}
 for all $t\geq 0$. 

Now, from \cite{Lei20}, there exist kernels $\rho(x, y)$ and $\tilde \psi(x)$ such that  the inverse transformation of (\ref{transwu}) has the form  
 \begin{align}\label{invtranswu}
   \hat{u}(x,t)&\!=\! \hat{w}(x,t)\!-\!\int_0^x \rho(x,y)\hat{w}(y,t)dy\!+\!\tilde\psi(x)\hat Z(t).
  \end{align}

\subsection*{Proof of Theorem 12(ii) - ISS with State and Error Estimation}
 Since transformations (\ref{transuw}) is invertible, and kernels $p$ and $q$ are bounded, then we have
 \begin{eqnarray}
     \|\tilde u\|\leq \tilde c \|\tilde w\|, && \|\tilde w\|\leq \tilde d \|\tilde u\|,
 \end{eqnarray}
 where $\tilde c$ and $\tilde d$ are two positive constants. Then, from (\ref{V21}) we obtain
 \begin{align}\label{V21}
\tilde\mu_1\left(|\tilde X|^2\!+\!\|\tilde u\|^2\right)\leq V_2\leq \tilde\mu_2\left(|\tilde X|^2\!+\!\|\tilde u\|^2\right),
\end{align}
for some positive constants $\tilde\mu_1$ and $\tilde\mu_2$.
From invertible transformation (\ref{transXZ})-(\ref{transwu}), we obtain
\begin{eqnarray}
\hat X(t)&=& \Delta_r^{-1} \hat Z(t),\label{transZX}\\
\hat u&=& \mathcal F^{-1}\left(\hat w+\psi(.)\hat Z\right),\label{transuw}
\end{eqnarray}
where $\mathcal F: L^2(0,1)\to L^2(0,1)$ is the invertible operator defined by 
$$\mathcal F(\hat u)(x,t)= \hat u(x, t)-\int_0^x s(x, y) \hat u(y,t)dy.$$
From (\ref{V1q}) and using the continuity of the operators $\mathcal F$ and $\mathcal F^{-1}$, and the boundedness of kernels $\psi$ and $s$, we establish 
\begin{equation}\label{V1qq}
        \hat\nu_1\left\|(\hat X, \hat u)\right\|_{\mathcal H}^2\leq V_1\leq \hat\nu_2\left\|(\hat X, \hat u)\right\|_{\mathcal H}^2,
        \end{equation}
where $\hat\nu_1$ and $\hat\nu_2$ are two positive constants. Putting all together gives 
\begin{align}
&c_1 \left\| \left( \hat{X}(t), \hat u(\cdot,t), \tilde{X}(t), \tilde u(\cdot,t) \right) \right\|_{\mathcal{H}^2}^2 \leq \mathcal{V}(t) \leq \notag\\
&c_2 \left\| \left( \hat{X}(t), \hat u(\cdot,t), \tilde{X}(t), \tilde u(\cdot,t) \right) \right\|_{\mathcal{H}^2}^2,\label{normequ}
\end{align}
for constants $c_1, c_2 > 0$. Thus, from the estimation (\ref{issclosed}), the norm equivalence (\ref{normequ}), and relations (\ref{errortrans1})-(\ref{errortrans2}) we obtain the ISS estimation (\ref{ISSthm}).
This completes the proof of Theorem 12.

  \end{pf}

\section{Conclusion}\label{sec5}

In this study, we have addressed the problem of global disturbance-to-state stabilization for nonlinear cascaded systems, consisting of ordinary differential equations coupled to a $1-D$ heat diffusion equation. 
 The control design method is based on a non-trivial combination of the high-gain observer, high-gain control, and backstepping method for coupled ODE-PDE systems introduced recently in \cite{Abd20} and \cite{AbdMoh21}. Establishing the existence and uniqueness of the solution to the closed-loop system posed an initial formidable challenge, which was surmounted by leveraging the powerful mathematical frameworks of $C_0$-semigroup theory, spectral decomposition, and the profound Bari theorem. Then, the DISS is obtained by Lypaunov analysis techniques. Building upon recent developments in \cite{dashkovskiy2023robust}, future research directions include extending these results to systems with Dirichlet boundary conditions and more general nonlinear terms.

 In recent times, there has been a significant amount of research conducted on the topic of finite time stabilization \cite{coron20} and delay systems \cite{prieur18}. The exploration of these areas has been driven by the recognition of their widespread practical application. Moving forward, our forthcoming work aims to integrate both of these concepts into our treated system.

\appendix

\section{Proof of Theorem \ref{Th1}}\label{pfTh1}

First, we  show that $\mathcal A+\lambda I$ is invertible and that $(\mathcal A+\lambda I)^{-1}$ is compact, for some positive $\lambda$.  

To do so,  let $[Y^T, v]^T\in \mathcal H$, we need to find a unique $[X^T, u]^T\in \mathcal D(\mathcal A)$ such
that $(\mathcal A+\lambda I)[X^T, u]^T = [Y^T, v]^T$. That is
\begin{eqnarray}
(A+\lambda I_n)X + B u(0) &=& Y,\label{XY}\\
u''+ (c+\lambda) u &=& v,\label{uv}\\
u'(0)=u'(1) &=& 0.\label{u0u1}
\end{eqnarray}
Using the theory of linear second order  ordinary differential equation, a general solution of   (\ref{uv}) and (\ref{u0u1})  exists and is expressed as follows
 \begin{eqnarray}\label{u}
  u(x)&=&\alpha \cos\left( \sqrt{ c\!+\!\lambda}x\right)\!+ \!\beta \cos\left( \sqrt{ c\!+\!\lambda}x\right)\!+\!u_p(x),
\end{eqnarray}
for all $ x\in [0, 1]$, where $\alpha$ and $\beta$ are real constants and $u_p(.)\in H^1(0, 1)$ (since $v\in L^2(0, 1)$) is a fixed particular solution of the forced system (\ref{uv}). If $\sqrt{c+\lambda}\neq k\pi$, for $k$ positive integer, then equations (\ref{uv}) and (\ref{u0u1}) admit a unique solution $u$ as (\ref{u}) and $\alpha, \beta$ are uniquely determined.
Then, since the function $u$ is unique and it is at least differentiable on $(0, 1)$, then $u(0)$ exists and it is unique.
Further, the matrix $A+\lambda I$  is invertible, meaning equation (\ref{XY}) has a unique solution $X=(A+\lambda I_n)^{-1}\left(Y-B u(0)\right)$, which implies $\mathcal A+\lambda I$ is invertible as well. According to Sobolev embedding theorem \cite{Adams} (Theorem 4.12, p. 85), since $\sigma(\mathcal A)$, the spectrum of  $\mathcal A$,  is nonempty, $(\mathcal A+\lambda I)^{-1}$  is compact on $\mathcal H$. So $\sigma(\mathcal A)=\sigma_p(\mathcal A)$, furthermore it  contains only  isolated eigenvalues \cite{Adams} (Theorem 2.16 page 27). Let us compute the point spectrum $\sigma_p(\mathcal A)$ of the operator $\mathcal A$. This turns to solve  the following  equations for $\lambda \in \sigma_p(\mathcal A)$
 \begin{eqnarray}
AY +B u(0)&=& \lambda Y,\label{eig1}\\
u''+ c u &=& \lambda u,\label{eig2}\\
u'(0)=u'(1) &=& 0.\label{eig3}
\end{eqnarray}
In the case when $u(x)=0, \forall x \in[0,1]$, the (\ref{eig1}) admits a unique solution $\lambda=0$. And $Y=\alpha e_1$, where $e_1=[1, 0, \cdots, 0]^T\in\mathbb R^n$ and $\alpha$ is real number. Then, $0\in \sigma_p(\mathcal A)$ and it is  associated to the eigenvector $\mathcal E_1=(e_1, 0)\in \mathcal H$.
Otherwise,  the unique solution of (\ref{eig2})-(\ref{eig3}) is given by
\begin{equation}
u(x)=u_k(x)=\alpha\cos(k\pi x), \; x\in[0, 1], 
\end{equation}
with $\; \lambda=\lambda_k=c-k^2\pi^2, k\in\mathbb N$, and $\alpha$ is a free constant. Then, (\ref{eig1}) gives $Y=\alpha Y_k$, where $Y_k=[1/\lambda_k^n, 1/\lambda_k^{n-1}, \cdots, 1/\lambda_k]^T$, where $k\in \mathbb N$. Thus, the spectrum $\sigma(\mathcal A)$ of the operator $\mathcal A$  is reduced to its point spectrum:  
$$\sigma(\mathcal A)=\sigma_p(\mathcal A)=\{0\}\cup\left\{\lambda_k=c-k^2\pi^2, \; k\in\mathbb N\right\},$$
whose eigenvectors are $\left\{\mathcal E_1, \mathcal Y_k=(Y_k, u_k), k\in \mathbb N\right\}.$

Let's show first that $\{e_1, Y_1, \cdots, Y_{n-1}\}$ form a basis of $\mathbb R^n$. Observe that 
\begin{eqnarray*}
\det(e_1, Y_1, \cdots, Y_{n-1})&=&\left|\begin{array}{ccccc}
1&1/\lambda_1^n & \cdots & 1/\lambda_{n-1}^{n-1}\\
0& 1/\lambda_1^{n-1} & \cdots & 1/\lambda_{n-1}^{n-1}\\
\vdots& \vdots &\ddots& \vdots\\
0& 1/\lambda_1 & \cdots& 1/\lambda_{n-1}^{n-1}
\end{array} 
\right|,\\
&=&\frac{1}{\prod_{i=1}^{n-1}\lambda_i}Vand(1/\lambda_1, \cdots, 1/\lambda_{n-1}),
\end{eqnarray*}
 where $Vand(x_1, \cdots, x_{n-1})$ is the Vandermonde determinant of $(x_1, \cdots, x_{n-1})^T\in\mathbb R^{n-1}$. Since $\lambda_i\neq \lambda_j$, for all $i\neq j$, it follows that $Vand(1/\lambda_1, \cdots, 1/\lambda_{n-1})\neq 0$ and then we conclude that $\det(e_1, Y_1, \cdots, Y_{n-1})\neq 0$ and this shows the desired result.

 Let $\mathcal U_k =(0_{\mathbb R^n}, u_k)\in \mathcal H, k\in \mathbb N$ and $\tilde{\mathcal Y}_k=(Y_k, 0) \in \mathcal H, 1\leq k\leq n-1$. Since $\left\{\mathcal U_k, k\in \mathbb N\right\}$ is a Riesz basis of $\{0_{\mathbb R^n}\}\times L^2(0, 1)$ and $\left\{\mathcal E_1, \tilde{\mathcal Y}_k, 1\leq k\leq n-1\right\}$ is a basis of the Euclidien space $\mathbb R^n\times \{0\}$, it follows that  $\left\{\mathcal E_1, \tilde{\mathcal Y}_1, \cdots, \tilde{\mathcal Y}_{n-1}, \mathcal U_k, k\in \mathbb N\right\}$ is Riesz basis of $\mathcal H$.

Now, Based on Bari Theorem \cite{Urata63}, we conclude that $\{\mathcal E_1, \mathcal Y_k, k\in \mathbb N\}$ forms a Riesz basis for $\mathcal H$. Indeed, we have 
\begin{eqnarray}
\|\mathcal E_1-\mathcal E_1\|^2+\sum_{i=1}^{n-1}\|\tilde{\mathcal Y}_i-\mathcal Y_i\|^2&+&\sum_{i=n}^{+\infty}\|\mathcal U_i-\mathcal Y_i\|^2\nonumber\\&=&\sum_{i=1}^{n-}\|u_i\|^2+\sum_{i=n}^{+\infty}|Y_i|^2,\nonumber\\
&=&\sum_{i=1}^{n-1}\|u_i\|^2+\sum_{i=n}^{+\infty}\sum_{j=1}^n\frac{1}{\lambda_i^{2j}},\nonumber\\
&\leq&\sum_{i=1}^{n-1}\|u_i\|^2+\tilde C\sum_{i=n}^{+\infty}\frac{1}{i^{4}},\nonumber\\
&<&+\infty,\label{bari}
\end{eqnarray}
where $\tilde C$ is a positive constant. Then, since $\{\mathcal E_1, \mathcal Y_k, k\in \mathbb N\}$ is $\omega$-linearly independent (as eigenvectors associated with different two-by-two eigenvalues), and $\left\{\mathcal E_1, \tilde{\mathcal Y}_1, \cdots, \tilde{\mathcal Y}_{n-1}, \mathcal U_k, k\in \mathbb N\right\}$ is Riesz basis of $\mathcal H$, using Bari Theorem, by (\ref{bari}) we obtain that $\{\mathcal E_1, \mathcal Y_k, k\in \mathbb N\}$ forms a Riesz basis for $\mathcal H$.

To summarise, the unbounded operator $\mathcal A$ has a discrete spectrum $\sigma(\mathcal A)=\sigma_p(\mathcal A)$, and its eigenvectors form a Riesz basis of $\mathcal H$.

Since all eigenvalues are bounded by above ($\lambda_k\leq c$, for all $k\in \mathbb N$), the operator $\mathcal A$ generates a $C_0$-semigroup $T(t)=e^{t\mathcal A}$ on $\mathcal H$. Furthermore, the nonlinear term $\mathcal F(\mathcal X)$ is locally Lipschitz, it is deduced that for all initial condition $\mathcal X_0\in \mathcal H$, there exists a unique mild solution of (\ref{dotX=AX++}) defined on a maximal interval $[0, T_{\max})$ by 
\begin{eqnarray}
\mathcal X(t)&=&T(t)\mathcal X_0+\int_0^t T(t-s)\mathcal F(\mathcal X(s))ds\notag\\&&+\int_0^tT(t-s)\mathcal  B_1U(s)ds \nonumber\\
&+&\!\!\!\int_0^t\!T(t\!-\!s)\!\!\left(\mathcal\!B_1d_4(s)\!+\! \mathcal B_2 \tilde d_1(s)\!+\!\mathcal B_3 d_3(s)\!\right)\!ds ,
\end{eqnarray}
In the following, we show that the operators $\mathcal B_1, \mathcal B_2$ and $\mathcal B_3$ are admissible for $T(t)$. To do so, we apply Proposition 3.9 in \cite{Mironchenko20}. First, we have 
\begin{eqnarray*}
    \int_0^tT(t-s)\mathcal B_1 U(s) ds&=&\int_0^te^{\lambda_{-1}(t-s)}<\mathcal B_1, \mathcal E_1>U(s)\mathcal E_1ds\\
    &+&\!\!\int_0^t\!\sum_{k=0}^{\infty}e^{\lambda_k(t-s)}\!\!< \mathcal B_1, \mathcal Y_k> U(s)\mathcal Y_kds,\\
&=&\sum_{k=0}^{\infty}u_k(1)\int_0^te^{\lambda_k(t-s)}U(s)ds\mathcal Y_k.
    \end{eqnarray*}
Then, we get
$$ \left\|\int_0^tT(t-s)\mathcal B_1 U(s) ds\right\|\leq N_t \int_0^tU^2(s)ds.$$
where $N_t\geq 0$. Moreover, we have
\begin{eqnarray*}
    \left\|\int_0^tT(t-s)\mathcal B_1 d_4(s) ds\right\|&\leq & H_t \sup_{s\in [0, t]}|d_4(s)|.
    \end{eqnarray*}
Furthermore, we have \begin{eqnarray*}
    \int_0^t\!\!\!T(t\!-\!s)\mathcal B_2 \tilde d_1(s)ds&=& \!\!\int_0^te^{\lambda_{-1}(t-s)}\!\!<\mathcal B_2 \tilde d_1(s), \mathcal E_1>ds \mathcal E_1\\
   &+&\!\!\int_0^t\!\!\sum_{k=0}^{\infty}e^{\lambda_k(t-s)}\!\!\!< \mathcal B_1, \mathcal Y_k\!>\! U(s)\mathcal Y_kds,\\
   &=& \int_0^te^{\lambda_{-1}(t-s)}e_1^TB_1 d_1(s) ds \mathcal E_1\\
   &+&\sum_{k=0}^{\infty}\int_0^te^{\lambda_k(t-s)}\Big( Y_k^TB_1 d_1(s)\notag\\
   &&\left.+\int_0^1d_2(x, s) u_k(x) dx \right)ds \mathcal Y_k.
    \end{eqnarray*}
Then, 
\begin{eqnarray*}
    \left\|\int_0^t\!\!T(t\!-\!s)\mathcal B_2 \tilde d_1(s) ds\right\|&\leq & K_t \sup_{s\in [0, t]}|d_1(s)|\\
    && +K_t\sup_{s\in [0, t]}\|d_2(., s)\|, 
    \end{eqnarray*}
    for some constant $K_t\geq 0$. The final term is given by
    \begin{eqnarray*}
    \int_0^t\!\!T(t\!-\!s)\mathcal B_3 \tilde d_3(s)ds&=&\!\!\int_0^t\!\!e^{\lambda_{-1}(t-s)}\!\!\!<\mathcal B_3 \tilde d_3(s),\! \mathcal E_1\!\!>\!ds \mathcal E_1\\
   &+&\!\!\int_0^t\!\!\sum_{k=0}^{\infty}\!\!\!e^{\lambda_k(t-s)}\!\!\!\!<\mathcal B_3 d_3(s),\!\mathcal Y_k\!\!> \!\mathcal Y_kds,\\
   &=& \sum_{k=0}^{\infty}\int_0^te^{\lambda_k(t-s)} d_3(s)ds \mathcal Y_k.
    \end{eqnarray*}
    Then, it follows
 \begin{eqnarray*}
    \left\|\int_0^tT(t-s)\mathcal B_3 \tilde d_3(s)ds\right\|&\leq & \sum_{k=0}^{\infty}\frac{1}{|\lambda_k|} \sup_{s\in [0, t]}|d_3(s)|\\
    &\leq& M  \sup_{s\in [0, t]}|d_3(s)|.
    \end{eqnarray*}   
Thus, operators $\mathcal B_1, \mathcal B_2,$ and $\mathcal B_3$ are admissible for $T(t)$. Since $\|\mathcal F(\mathcal X)\|\leq c \|\mathcal X\|$, for some positive constant and using Gronwall inequality, we obtain (\ref{priorest}),  for all $t\in [0, T_{\max})$, where $c_t$ and $d_t$ are positive reals. Thus, the proof of Theorem \ref{Th1} is established.

\section{ Proof o Lemma \ref{invtransuw} } \label{proofinvtransuw}   From (\ref{transuw}), we have
\begin{align}
\tilde u_t(x, t)&=\tilde w_t(x, t)-\int_x^1q(x, y) \tilde w_t(y, t)dy\notag\\
&=\tilde w_{xx}(x,t)+\tilde k(x)\tilde w(1, t)+\tilde d_2(x, t)\notag\\
&-  \int_x^1q(x, y)\big( \tilde w_{xx}(y, t)+\tilde k(y)\tilde w(1, t)\notag\\
&+ \tilde d_2(y, t)\big)dy.\label{ut}
\end{align}
Twice integration by parts, it follows
\begin{align}
 \int_x^1\!\!q(x,\! y) \tilde w_{xx}(y,\!t)dy&=q(x,\!1) \tilde w_{x}(1,\!t)\!-\!q(x,\!x)\!\tilde w_{x}(x,\!t)\notag\\
 &-q_y(x, 1) \tilde w(1,t)+q_y(x, x) \tilde w(x,t)\notag\\
 &+\int_x^1q_{yy}(x, y) \tilde w(y,t)dy.\label{intpart}
\end{align}
Putting (\ref{intpart}) in (\ref{ut}), it yields
\begin{align}
\tilde u_t(x, t)&=\tilde w_{xx}(x,t)+\tilde k(x)\tilde w(1, t)+\tilde d_2(x, t)\notag\\
&-q(x, 1) \tilde w_{x}(1,t)+q(x, x) \tilde w_{x}(x,t)\notag\\&+q_y(x, 1) \tilde w(1,t)-q_y(x, x) \tilde w(x,t)\notag\\
&-\int_x^1q_{yy}(x, y) \tilde w(y,t)dy\notag\\
&-\int_x^1q(x, y)\tilde k(y)dy \; \tilde w(1, t)\notag\\
&-\int_x^1q(x, y) \tilde d_2(y, t)dy.
\end{align} 
On the other hand, from (42), we get
\begin{align}
\tilde u_x(x, t)&=\tilde w_{x}(x,t)+q(x, x) \tilde w(x, t)\notag\\
&-\int_x^1q_{x}(x, y) \tilde w(y, t)dy,\notag\\
\tilde u_{xx}(x, t)&=\tilde w_{xx}(x,t)+(q(x, x))' \tilde w(x, t)\notag\\
&+q(x, x) \tilde w_x(x, t)+q_x(x, x) \tilde w(x, t)\notag\\
&-\int_x^1q_{xx}(x, y) \tilde w(y, t)dy.
\end{align} 
Then, it yields
\begin{align}
\tilde u_t-\tilde u_{xx}-c\tilde u&=\tilde w_{xx}+\tilde k(x)\tilde w(1, t)+ \tilde d_2(x, t)\notag\\
&- q(x, 1) \tilde w_{x}(1,t)+q(x, x) \tilde w_{x}\notag\\
&+q_y(x, 1) \tilde w(1,t)-q_y(x, x) \tilde w\notag\\
&-\int_x^1q_{yy}(x, y) \tilde w(y,t)dy\notag\\
&-\int_x^1q(x, y)\tilde k(y)dy \; \tilde w(1, t)\notag\\
&-\int_x^1q(x, y) \tilde d_2(y, t)dy-\tilde w_{xx}\notag\\
&-(q(x, x))' \tilde w-q(x, x) \tilde w_x\notag\\
&-q_x(x, x) \tilde w-c\tilde w\notag\\&+\int_x^1q_{xx}(x, y) \tilde w(y, t)dy\notag\\&+c\int_x^1q(x, y)\tilde w(y, t)dy,\notag\\
&= \big(-c-2(q(x,x))'\big)\tilde w(x,t)\notag\\
&+\!\Big(\tilde k(x)\!+\!q_y(x, \!1)\!-\!\int_x^1\!\!q(x,\!y)\tilde k(y)dy \Big)\tilde w(1,\!t)\notag\\
&+\int_x^1\left(q_{xx}-q_{yy}+cq\right)\tilde w(y, t)dy\notag\\
&-\int_x^1q(x, y)\tilde d_2(y, t)dy\notag\\
&+\tilde d_2(x, t)-q(x, 1) \tilde w_x(1, t).\label{utttt}
\end{align}
Furthermore, we have
\begin{eqnarray}
\tilde w_x(1, t)&=&\tilde u_x(1, t)-q(1, 1)\tilde w(1, t),\notag\\
&=&\tilde u_x(1, t)-q(1, 1)\tilde u(1, t),\notag\\
&=&-(q_2+q(1, 1))\tilde u(1, t)+d_4(t).\label{wxxx}
\end{eqnarray}
In view of (\ref{pxx1}), (\ref{pxx2}) and (\ref{wxxx}), from (\ref{utttt}), it follows that
\begin{align*}
\tilde u_t-\tilde u_{xx}-c\tilde u&= \Big( q(x, 1)(q_2+q(1, 1))+\tilde k(x)\notag\\
&+q_y(x, 1)-\int_x^1q(x, y)\tilde k(y)dy \Big)\tilde u(1, t)\\
&+\tilde d_2(x, t)-q(x, 1)d_4(t)\notag\\
&-\int_x^1q(x, y)\tilde d_2(y, t)dy.
\end{align*}
To guarantee that $\tilde u(x, t)$ satisfies equation (39), our task is simply to identify values for $\tilde k(\cdot)$ and $\tilde d_2(\cdot, \cdot)$ that fulfill the subsequent equations
\begin{align}\label{kkk}
    -k(x)&= q(x, 1)(q_2+q(1, 1))+\tilde k(x)+q_y(x, 1)\notag\\
    &-\int_x^1q(x, y)\tilde k(y)dy \\ \label{ddd}
d_2(x, t)&=\tilde d_2(x, t)-q(x, 1)d_4(t)\notag\\
&-\int_x^1q(x, y)\tilde d_2(y, t)dy.
\end{align}
To accomplish this, it only need to show the surjectivity of the following integral operator  
$$\mathcal Q(\beta)(x,t)= \beta(x, t)-\int_x^1q(x, y) \beta(y, t) dy, $$
for all $x\in [0, 1]$ and $t\geq 0$.  To do this, for $\gamma(\cdot, \cdot)\in C([0, 1]\times [0, +\infty))$, we must find $\beta(\cdot, \cdot)\in C([0, 1]\times [0, +\infty))$ such that $\mathcal Q(\beta)=\gamma$, which is equivalent to
\begin{equation}
    \gamma(x, t)=\beta(x, t)-\int_x^1q(x, y) \beta(y, t) dy,
\end{equation}
for $(x, t)\in [0, 1]\times [0, +\infty)$. This is can be established  using  the classical Picard method of successive approximation. In fact, let $\beta_0(x, t)=0$ and 
$$\beta_{n+1}(x, t)=\gamma(x, t)+\int_x^1q(x, y) \beta_n(y, t) dy, \;\forall x\in [0, 1], \forall t\geq 0.$$
Define $\Delta \beta_n(x, t)=\beta_{n+1}(x, t)-\beta_n(x, t)$. We obtain $\Delta \beta_0(x, t)=\gamma(x, t)$ and 
$$\Delta\beta_{n+1}(x, t)=\int_x^1q(x, y)\Delta\beta_n(y, t) dy.$$
It can be proved by simple recursive method that, for all $n\in \mathbb N$,
we have
$$\left|\Delta\beta_n(x, t)\right|\leq M(t)^{n+1}\frac{(1-x)^n}{n!}, \; \forall x\in [0, 1],$$
where $M(t)=\sup_{(x, y, s)\in \mathcal T_2\times [0, t]}\left| q(x, y)\gamma(y, s)\right|$. It then follows that the series 
$$\sum_{k=0}^{+\infty}\Delta\beta_k(x, t)=\lim_{n\rightarrow+\infty}\beta_n(x, t)=\beta(x, t)$$
converges uniformly in any compact set included in $[0, 1]\times [0, +\infty)$. Then, the surjectivity of the operator $\mathcal Q$ is proved. Thus, the existence of $k(\cdot)$ and $\tilde d_2(\cdot, \cdot)$ satisfying (\ref{kkk}), and (\ref{ddd}) is done, and the proof of Lemma \ref{invtransuw} is finished.

\section{Proof of Lemma \ref{lemk}}\label{pflemk}
Let solve the following integral equation using the standard Picard method of successive approximation 
$$k(x)= h(x)-\int_x^1p(x, y) k(y) dy, $$
for all $x\in [0, 1]$, where $h(x)$ is a given continuous function and $k(x)$ is the unknown function to be determined.    In fact, let $k_0(x)=0$ and 
$$k_{n+1}(x)=h(x)+\int_x^1p(x, y) k_n(y) dy, \;\forall x\in [0, 1].$$
Define $\Delta k_n(x)=k_{n+1}(x)-k_n(x)$. We obtain $\Delta k_0(x)=h(x)$ and 
$$\Delta k_{n+1}(x)=\int_x^1p(x, y)\Delta k_n(y) dy.$$
It can be proved by simple recursive method that, for all $n\in \mathbb N$,
we have
$$\left|\Delta k_n(x)\right|\leq M^{n+1}\frac{(1-x)^n}{n!}, \; \forall x\in [0, 1],$$
where $M=\sup_{(x, y)\in \mathcal T_2}\left| p(x, y)\gamma(y, s)\right|$. It then follows that the series 
$$\sum_{n=0}^{+\infty}\Delta k_n(x)=\lim_{n\rightarrow+\infty} k_n(x)=k(x)$$
converges uniformly in $[0, 1]$. Then, the integral equation (\ref{k}) has a continuous solution $k(x), x\in [0, 1]$.

\section{Proof of Lemma \ref{V1}} \label{lemm3pf} 
Concerning the first point, it is straightforward to establish (\ref{V1q}) with $\nu_1=\min(\lambda_{\min}(P_1), a/2)$ and $\nu_2=\max(\lambda_{\max}(P_1), a/2)$, where $\lambda_{\min}(P_1)$ and $\lambda_{\max}(P_1)$ are the smallest and the largest eigenvalues of the matrix $P_1$ respectively. Now let's proceed to the proof of the second item of the Lemma \ref{V1}. The derivative of $V_1(t)$ along the solutions of system (\ref{closeddimfin1})-(\ref{closedbound22}) satisfies
\begin{eqnarray}\label{dotV1}
\dot V_1(t)&\leq & -r^{-1}\delta_1|\hat Z(t)|^2+2r^{-1}\hat Z(t)^TP_1B\hat w(0, t)\notag\\
&-&\!2r^{-1}\hat Z^T(t)P_1LC\tilde X(t)-a\|\hat w_x(., t)\|^2\!-\!a\eta\hat w(1, t)^2 \nonumber\\
&+&a r^{-1}\int_0^1\hat w(y, t)\psi(y)dy LC\tilde X(t).
\end{eqnarray}
Let's adequately bound the cross terms in (\ref{dotV1}). Completing squares, we get
\begin{eqnarray}
2r^{-1}\hat Z(t)^TP_1B\hat w(0, t)&\leq & r^{-2}|P_1|^2|\hat Z(t)|^2\!+\! \hat w(0, t)^2,\label{eq1}\\
-2\hat Z^T(t)P_1LC\tilde X(t)&\leq &|\hat Z(t)|^2\!+\!  |P_1L|^2 |\tilde X(t)|^2,\label{eq2}
\end{eqnarray}
where we have used that $|B|=|C|=1$.
In the same spirit, and using (\ref{cpsi}), it follows
\begin{flalign}
\int_0^1\!\!\!\hat w(y, t)\psi(y)dy LC\tilde X(t)\leq &\frac{c_{\psi}|L|}{2} \!\! \left( \|\hat w(., t)\|^2\!+\!|\tilde X(t)|^2\!\right).\label{eq6}
\end{flalign}
Using the fact that $\hat w(0, t) = \hat w(1, t)-\int^1_0 \hat w_x(x, t)dx$ and by applying Cauchy-Schwarz inequality, it follows
\begin{equation}\label{hatv}
\hat w(0, t)^2 \leq 2\hat w(1, t)^2+2\|\hat w_x(\cdot, t)\|^2,
\end{equation} 
Substituting inequalities (\ref{eq1})-(\ref{eq6}) and (\ref{hatv}) in (\ref{dotV1}), it follows that
\begin{eqnarray}
\dot V_1(t)&\leq&\!-\!r^{-1}\!\!\left(\delta_1\!-\!r^{-1}\!|P_1|^2\!-\!1\right)|\hat Z(t)|^2\!-\!\left(a\!-\!2\right)\|\hat w_x(\cdot,\!t)\|^2\nonumber\\
&&+\frac{a}{2}r^{-1}c_{\psi}|L|\|\hat w(., t)\|^2-(a\eta-2)\hat w(1, t)^2\nonumber\\
&& + \left(r^{-1}|P_1L|^2+\frac{1}{2}ar^{-1}c_{\psi} |L|\right) |\tilde X(t)|^2.\label{dotV11}
\end{eqnarray}
Now, for $a>2$, by  using the inequality (which is derived from Poincarré inequality) $$-(a-2)\|\hat w_x(., t)\|^2\leq \frac{1}{2}(a-2)\hat w(1, t)^2-\frac{1}{4}(a-2)\|\hat w(., t)\|^2,$$ 
inequality (\ref{dotV11}) becomes (\ref{dotV112}).

\section{Proof of Lemma \ref{V2}}\label{lemma5prf}
The derivative of $V_2(t)$ along the solutions of system (\ref{closeddimfin2}), (\ref{closedheat2}), (\ref{closedbound2}), (\ref{closedbound22}) satisfies
\begin{eqnarray}\label{dotV21}
\dot V_2(t)&\leq & -r^{-1}\delta_2 |\tilde X(t)|^2+2r^{-1}\tilde X(t)^TP_2B\tilde w(0, t)\notag\\
&&-2r^{-1}\tilde X^T(t)P_2B\int_0^1q(0, y)\tilde w(y, t)dy,\nonumber\\
&&+2\tilde X^T(t)P_2\tilde f(X(t))+2\tilde X(t)^TP_2 \Delta_rB_1d_1(t)\notag\\
&&-b\|\tilde w_x(., t)\|^2-b(q_2-c/2)\tilde w(1, t)^2 \nonumber\\
&&-b \tilde w(0, t) d_3(t)+b\tilde w(1, t)d_4(t)\nonumber\\
&&+b \int_0^1\tilde w(y, t)\tilde d_2(y, t) dy.
\end{eqnarray}
In view of Assumption \ref{ass1} and the fact that the scalar gain $r\geq 1$, the $i$-th component $\tilde f_i(X)$ of the nonlinear term $\tilde f(X(t))=\Delta_r f(\Delta_r^{-1}(\hat Z(t)+\tilde X(t)))$  satisfies 
\begin{eqnarray}
\left|\tilde f_i(X(t))\right|&\leq &\theta r^{-n+i-1}\sum_{j=i+2}^nr^{n-j+1}\left|\hat Z_j+\tilde X_j \right|,\nonumber\\
&\leq &\theta \sum_{j=i+2}^nr^{i-j}\left|\hat Z_j+\tilde X_j \right|,\nonumber\\
&\leq & \theta r^{-2}\sum_{j=i+2}^n\left(\left|\hat Z_j\right|+\left|\tilde X_j \right|\right), \; (\textsf{since } r\geq 1),
\nonumber\\
&\leq & \theta r^{-2}\sqrt{n\!-\!i\!+\!1}\left(\!\!\left(\sum_{j=i+2}^n\!\hat Z_j^2\right)^{\frac{1}{2}}\!\!\!+\!\!\!\left(\sum_{j=i+2}^n\tilde X_j^2\!\!\right)^{\frac{1}{2}}\!\!\right),\nonumber\\
&\leq & \theta r^{-2}\sqrt n\left(\left|\hat Z\right|+\left|\tilde X\right|\right).\label{fi}
\end{eqnarray}
Using (\ref{fi}), it yields
\begin{eqnarray}
\left|\tilde f(X(t))\right|&=&\sqrt{\sum_{i=1}^n\tilde f_i(X(t))^2},\nonumber\\
&\leq & \sqrt{\theta^2r^{-4} n\sum_{i=1}^n\left(\left|\hat Z(t)\right|+\left|\tilde X(t)\right|\right)^2},\nonumber\\
&\leq & n \theta r^{-2}\left(\left|\hat Z(t)\right|+\left|\tilde X(t)\right|\right).\label{fifi}
\end{eqnarray}
Then, by (\ref{fifi}), we obtain
\begin{align}
    &2r^{-1}\!\!\tilde X(t)^T\!P_2\!B\tilde w(0, t)\leq\!\! \lambda^2r^{-2}|\tilde X(t)|^2\!+\!\tilde w(0, t)^2,\label{a1}\\
&2\tilde X^T(t)P_2\tilde f(X(t))\leq \!\!\frac{\lambda n\theta}{r^2} \left(3|\tilde X(t)|^2\!\!+\!\!|\hat Z(t)|^2\right),\label{a2}\\
&-2r^{-1}\tilde X^T(t)P_2\!\!\int_0^1\!\!q(0,\!y)\tilde w(y, t)dy \leq \notag\\
& \hspace{2.5cm}\lambda^2r^{-2}|\tilde X(t)|^2+m_1\|\tilde w(., t)\|^2\label{a3},\\
&2\tilde X(t)^TP_2\Delta_r B_1 d_1(t)\leq \frac{\lambda^2}{r^2}|B_1|^2|\tilde X(t)|^2\!+\!d_1(t)^2,\\
&b\tilde w(1, t)d_4(t)\leq  \frac{1}{2}b^2\tilde w(1, t)^2+\frac{1}{2}d_4(t)^2,\\
&2b\!\!\int_0^1\!\!\tilde w(y, \!t)\tilde d_2(y,\!t)dy\leq \|\tilde w(.,\!t)\|^2\!+\!b^2 \|\tilde d_2(., t)\|^2. \label{a4}
\\
&-b\tilde w(0, t) d_3(t)\leq \frac{1}{2}\tilde w(0, t)^2+\frac{1}{2}b^2 d_3(t)^2,\label{a5}
\end{align}
where $\lambda=\lambda_{\max}(P_2).$
Putting (\ref{a1})-(\ref{a5}) in (\ref{dotV21}), it follows
\begin{align}\label{dotV22}
\dot V_2(t)&\leq  -r^{-1}\left(\delta_2- r^{-1} m_2\right) |\tilde X(t)|^2+\frac{3}{2}\tilde w(0, t)^2\notag\\
&+(m_1+\frac{1}{2})\|\tilde w(., t)\|^2-b\|\tilde w_x(., t)\|^2\notag\\
&-\left(b(q_2-c/2)-\frac{1}{2}b^2\right)\tilde w(1, t)^2\notag\\
&+ n\theta\lambda_{\max}(P_2)r^{-2}|\hat Z(t)|^2+d_1^2(t)\nonumber\\
& +\frac{1}{2}d_4(t)^2+\frac{1}{2}b^2\|\tilde d(., t)\|^2+\frac{1}{2}b^2 d_3(t)\|^2.
\end{align}
It is clear that $\tilde w(0, t) =\tilde w(1, t)-\int^1_0\tilde w_x(x, t)dx$, then the following useful
inequalities can be obtained by Schwarz inequality:
\begin{equation}\label{wwx}
\frac{3}{2}\tilde w(0, t)^2\leq 3\tilde w(1, t)^2+3\|\tilde w_x(., t)\|^2.
\end{equation}
Then, we get

\begin{align}\label{dotV23}
\dot V_2(t)&\leq  
-r^{-1}\!\left(\delta_2\!- \!r^{-1} m_2\right) |\tilde X(t)|^2\!+\!(m_1\!+\!\frac{1}{2})\|\tilde w(., t)\|^2\nonumber\\
&-(b\!-\!3)\|\tilde w_x(.,\!t)\|^2\!-\!\left(b(q_2\!-\!\frac{c}{2})-\frac{b^2}{2}\!-\!3\right)\tilde w(1,\!t)^2\nonumber\\
& + n\theta\lambda r^{-2}|\hat Z(t)|^2+d_1^2(t)+\frac{1}{2} b^2d_3(t)^2+\frac{1}{2}d_4(t)^2\notag\\
&+\frac{1}{2}b^2\|\tilde d(., t)\|^2,
\end{align}
Furthermore, based on the following inequality 
\begin{align}
  -\!(b\!-\!3)\|\tilde  w_x(.,\! t)\|^2\!\leq \!\frac{b\!-\!3}{2}\tilde w(1,\!t)^2\!-\!\frac{b\!-\!3}{4}\|\tilde w(.,\!t)\|^2,
  \end{align} 
we obtain (\ref{dotV2121}).




\begin{thebibliography}{99}     %

\bibitem{Adams}
Robert~A. Adams and John J.~F. Fournier.
\newblock {\em Sobolev Spaces}, volume 140 of {\em Pure and Applied
  Mathematics}.
\newblock Elseiver, New York, 2 edition, 2003.

\bibitem{Sontag08}
Andrei~A Agrachev, A~Stephen Morse, Eduardo~D Sontag, H{\'e}ctor~J Sussmann,
  Vadim~I Utkin, and Eduardo~D Sontag.
\newblock Input to state stability: Basic concepts and results.
\newblock {\em Nonlinear and optimal control theory: lectures given at the CIME
  summer school held in Cetraro, Italy June 19--29, 2004}, pages 163--220,
  2008.

\bibitem{Ahmed15}
T.~Ahmed-Ali, F.~Giri, M.~Krstic, and F.~Lamnabhi-Lagarrigue.
\newblock Observer design for a class of nonlinear ode–pde cascade systems.
\newblock {\em Systems \& Control Letters}, 83:19--27, 2015.

\bibitem{Praly09}
V.~Andrieu and L.~Praly.
\newblock A unifying point of view on output feedback designs for global
  asymptotic stabilization.
\newblock {\em Automatica}, 45(8):1789--1798, 2009.

\bibitem{Bastin21}
Georges Bastin, Jean-Michel Coron, and Amaury Hayat.
\newblock Input-to-state stability in sup norms for hyperbolic systems with
  boundary disturbances.
\newblock {\em Nonlinear Analysis}, 208:112300, 2021.

\bibitem{Abd20}
Abdallah Benabdallah.
\newblock Stabilization of a class of nonlinear uncertain ordinary differential
  equation by parabolic partial differential equation controller.
\newblock {\em Int. J. Robust Nonlinear Control}, 30:3023–3038, 2020.

\bibitem{AbdMoh21}
Abdallah Benabdallah and Mohsen Dlala.
\newblock Rapid exponential stabilization by boundary state feedback for a
  class of coupled nonlinear ode and $ 1-d $ heat diffusion equation.
\newblock {\em Discrete and Continuous Dynamical Systems - S},
  15(5):1085--1102, 2022.

\bibitem{Boyd94}
Stephen Boyd, Laurent El~Ghaoui, Eric Feron, and Venkataramanan Balakrishnan.
\newblock {\em Linear Matrix Inequalities in System and Control Theory}.
\newblock Society for Industrial and Applied Mathematics, 1994.

\bibitem{Chang21}
Yanjie Chang, Tongjun Sun, Xianfu Zhang, and Xiandong Chen.
\newblock Output feedback stabilization for a class of cascade nonlinear
  ode-pde systems.
\newblock {\em International Journal of Control, Automation and Systems},
  19(7):2519--2528, 2021.

\bibitem{coron20}
Jean-Michel Coron and Hoai-Minh Nguyen.
\newblock Finite-time stabilization in optimal time of homogeneous quasilinear
  hyperbolic systems in one dimensional space.
\newblock {\em ESAIM: Control, Optimisation and Calculus of Variations},
  26:119, 2020.

\bibitem{dashkovskiy2023robust}
Sergey Dashkovskiy, Oleksiy Kapustyan, and Vitalii Slynko.
\newblock Robust stability of a nonlinear ode-pde system.
\newblock {\em SIAM Journal on Control and Optimization}, 61(3):1760--1777,
  2023.

\bibitem{Mohsen22}
Mohsen Dlala and Abdallah Benabdallah.
\newblock Global stabilization of nonlinear finite dimensional system with
  dynamic controller governed by $1-d$ heat equation with neumann
  interconnection.
\newblock {\em Mathematics}, 10(2), 2022.

\bibitem{Jacob22}
R.~Hosfeld, B.~Jacob, and F.L. Schwenninger.
\newblock Integral input-to-state stability of unbounded bilinear control
  systems.
\newblock {\em Mathematics of Control, Signals, and Systems}, 34:273–295,
  2022.

\bibitem{Orliz18}
Birgit Jacob, Robert Nabiullin, Jonathan~R. Partington, and Felix~L.
  Schwenninger.
\newblock Infinite-dimensional input-to-state stability and orlicz spaces.
\newblock {\em SIAM Journal on Control and Optimization}, 56(2):868--889, 2018.

\bibitem{Jacoub19}
Birgit Jacob, Felix~L. Schwenninger, and Hans Zwart.
\newblock On continuity of solutions for parabolic control systems and
  input-to-state stability.
\newblock {\em Journal of Differential Equations}, 266(10):6284--6306, 2019.

\bibitem{Karf18}
Iasson Karafyllis and Miroslav Krstic.
\newblock {\em Input-to-State Stability for PDEs}.
\newblock Communications and Control Engi- neering, Springer, 2018.

\bibitem{Karafyllis19}
Iasson Karafyllis and Miroslav Krstic.
\newblock {\em Input-to-State Stability for PDEs. Communications and Control
  Engineering}.
\newblock Springer, 2019.

\bibitem{Kris09heat}
Miroslav Krstic.
\newblock Compensating actuator and sensor dynamics governed by diffusion pdes.
\newblock {\em Systems \& Control Letters}, 58(5):372--377, 2009.

\bibitem{Hao07}
Hao Lei and Wei Lin.
\newblock Adaptive regulation of uncertain nonlinear systems by output
  feedback: A universal control approach.
\newblock {\em Systems \& Control Letters}, 56(7):529--537, 2007.

\bibitem{Lei20}
Yuan Lei, Xinglan Liu, and Chengkang Xie.
\newblock Stabilization of an ode-pde cascaded system by boundary control.
\newblock {\em Journal of the Franklin Institute}, 357(14):9248--9267, 2020.

\bibitem{Lhachemi22}
Hugo Lhachemi and Christophe Prieur.
\newblock Stability analysis of reaction–diffusion pdes coupled at the
  boundaries with an ode.
\newblock {\em Automatica}, 144:110465, 2022.

\bibitem{Lhachemi20}
Hugo Lhachemi, David Saussié, Guchuan Zhu, and Robert Shorten.
\newblock Input-to-state stability of a clamped-free damped string in the
  presence of distributed and boundary disturbances.
\newblock {\em IEEE Transactions on Automatic Control}, 65(3):1248--1255, 2020.

\bibitem{Xie20}
Xinglan Liu and Chengkang Xie.
\newblock Control law in analytic expression of a system coupled by
  reaction–diffusion equation.
\newblock {\em Systems $\&$ Control Letters}, 137:104643, 2020.

\bibitem{Mazenc1996adding}
F.~Mazenc and L.~Praly.
\newblock Adding integrations, saturated controls, and stabilization for
  feedforward systems.
\newblock {\em IEEE Transactions on Automatic Control}, 41(11):1559--1578,
  1996.

\bibitem{Karafilis19}
A.~Mironchenko, I.~Karafyllis, and M.~Krstic.
\newblock Monotonicity methods for input-to-state stability of nonlinear
  parabolic pdes with boundary disturbances.
\newblock {\em SIAM J Control Optim}, 57:510–532, 2019.

\bibitem{Prieur20}
A.~Mironchenko and C.~Prieur.
\newblock Input-to-state stability of infinite-dimensional systems: recent
  results and open questions.
\newblock {\em SIAM Review}, 62(3):529--614, 2020.

\bibitem{Mironchenko19}
Andrii Mironchenko, Iasson Karafyllis, and Miroslav Krstic.
\newblock Monotonicity methods for input-to-state stability of nonlinear
  parabolic pdes with boundary disturbances.
\newblock {\em SIAM Journal on Control and Optimization}, 57(1):510--532, 2019.

\bibitem{Mironchenko20}
Andrii Mironchenko and Christophe Prieur.
\newblock Input-to-state stability of infinite-dimensional systems: Recent
  results and open questions.
\newblock {\em SIAM Review}, 62(3):529--614, 2020.

\bibitem{Pazy}
A.~Pazy.
\newblock {\em Semigroups of Linear Operators and Applications to Partial
  Differential Equations}.
\newblock Springer-Verlag, New York, 1992.

\bibitem{prieur18}
Christophe Prieur and Emmanuel Tr{\'e}lat.
\newblock Feedback stabilization of a 1-d linear reaction--diffusion equation
  with delay boundary control.
\newblock {\em IEEE Transactions on Automatic Control}, 64(4):1415--1425, 2018.

\bibitem{Sch20}
F.L. Schwenninger.
\newblock Input-to-state stability for parabolic boundary control:linear and
  semilinear systems.
\newblock In J.~Kerner, H.~Laasri, and D.~Mugnolo, editors, {\em Control Theory
  of Infinite-Dimensional Systems}, volume 277, pages 83--116. Birkhäuser,
  Cham., Berlin, 2020.

\bibitem{Kris04}
Krstic~M. Smyshlyaev~A.
\newblock Closed-form boundary state feedbacks for a class of $1-d$ partial
  integro-differential equations.
\newblock {\em IEEE Transaction on Automatic Control}, 49(12):2185--2202, 2004.

\bibitem{Tang11}
Shuxia Tang and Chengkang Xie.
\newblock Stabilization for a coupled pde–ode control system.
\newblock {\em Journal of the Franklin Institute}, 348(8):2142--2155, 2011.

\bibitem{Tang11SCL}
Shuxia Tang and Chengkang Xie.
\newblock State and output feedback boundary control for a coupled pde–ode
  system.
\newblock {\em Systems \& Control Letters}, 60(8):540--545, 2011.

\bibitem{Urata63}
Yosihiro Urata.
\newblock {A theorem of Bari on the completeness of orthonormal systems}.
\newblock {\em Proceedings of the Japan Academy}, 39(3):160--161, 1963.

\bibitem{Wang21}
Jun-Min Wang, Han-Wen Zhang, and Xiu-Fang Yu.
\newblock Input-to-state stabilization of coupled parabolic pdes subject to
  external disturbances.
\newblock {\em IMA Journal of Mathematical Control and Information},
  39(1):185--218, 12 2021.

\bibitem{Wonham67}
W.~Wonham.
\newblock On pole assignment in multi-input controllable linear systems.
\newblock {\em IEEE Transactions on Automatic Control}, 12(6):660--665, 1967.

\bibitem{Zhang21}
Han-Wen Zhang, Jun-Min Wang, and Jian-Jun Gu.
\newblock Exponential input-to-state stabilization of an ode cascaded with a
  reaction–diffusion equation subject to disturbances.
\newblock {\em Automatica}, 133:109885, 2021.

\bibitem{ZhangCDC21}
Hanwen Zhang, Junmin Wang, and Jianjun Gu.
\newblock Input-to-state stabilization for an ode cascaded by a parabolic pide
  with disturbances.
\newblock In {\em 2021 33rd Chinese Control and Decision Conference (CCDC)},
  pages 979--984, 2021.

\bibitem{Zhang201}
Yu-Long Zhang, Jun-Min Wang, and Donghai Li.
\newblock Input-to-state stabilization of an ode-wave system with disturbances.
\newblock {\em Mathematics of Control, Signals, and Systems}, 32(4):489--515,
  2020.

\bibitem{Zhang20}
Yu-Long Zhang, Jun-Min Wang, and Donghai Li.
\newblock Input-to-state stabilization of an ode-wave system with disturbances.
\newblock {\em Mathematics of Control, Signals, and Systems}, 32(4):489--515,
  2020.

\bibitem{Xie21}
Zhiyuan Zhen, Yuanchao Si, and Chengkang Xie.
\newblock Indirect control to stabilise reaction–diffusion equation.
\newblock {\em International Journal of Control}, 94(11):3091--3098, 2021.

\bibitem{zheng2018iss}
Jun Zheng, Hugo Lhachemi, Guchuan Zhu, and David Saussi{\'e}.
\newblock Iss with respect to boundary and in-domain disturbances for a coupled
  beam-string system.
\newblock {\em Mathematics of Control, Signals, and Systems}, 30(4):21, 2018.

\bibitem{Zheng18}
Jun Zheng and Guchuan Zhu.
\newblock Input-to-state stability with respect to boundary disturbances for a
  class of semi-linear parabolic equations.
\newblock {\em Automatica}, 97:271--277, 2018.

\bibitem{Zheng20}
Jun Zheng and Guchuan Zhu.
\newblock Input-to-state stability for a class of one-dimensional nonlinear
  parabolic pdes with nonlinear boundary conditions.
\newblock {\em SIAM Journal on Control and Optimization}, 58(4):2567--2587,
  2020.

\end{thebibliography}

\end{document}